\documentclass[11pt]{amsart}
\usepackage{latexsym}

\DeclareMathOperator{\ad}{ad}

\newcommand{\F}{\mathbb F}

\newcommand{\Z}[0]{\mathbb Z}
\newcommand{\q}{\bar{q}}

\newtheorem{dummy}{Dummy}

\numberwithin{dummy}{section}
\numberwithin{equation}{section}

\newtheorem{lemma}[dummy]{Lemma}

\newtheorem{theorem}[dummy]{Theorem}

\newtheorem{prop}[dummy]{Proposition}
\newtheorem{cor}[dummy]{Corollary}

\theoremstyle{definition}

\theoremstyle{remark}
\newtheorem{rem}[dummy]{Remark}

\begin{document}

\bibliographystyle{amsalpha}
\author{Marina Avitabile}
\email{marina.avitabile@unimib.it}
\address{Dipartimento di Matematica e Applicazioni\\
  Universit\`a di Milano-Bicocca\\
  via Cozzi 53\\
  I-20125 Milano\\
  Italy}

\author{Giuseppe Jurman}
\email{jurman@fbk.eu}
\address{FBK-irst\\
  via Sommarive 18\\
  I-38100 Povo (Trento)\\
  Italy}

\author{Sandro Mattarei}
\email{mattarei@science.unitn.it}
\address{Dipartimento di Matematica\\
  Universit\`a degli Studi di Trento\\
  via Sommarive 14\\
  I-38100 Povo (Trento)\\
  Italy}

\thanks{The first and third authors are members of INdAM-GNSAGA, Italy.
  They acknowledge financial  support from Ministero dell'Istruzione,
  dell'Universit\`a e della  Ricerca, Italy, to  the project
  ``Lie rings and algebras, groups, cryptography''.}

\title[The initial structure of thin Lie algebras]{The structure of thin {L}ie algebras\\ up to the second diamond}
\date{\today}
\begin{abstract}
{\em Thin} Lie algebras are graded Lie algebras $L=\bigoplus_{i=1}^{\infty}L_{i}$
with $\dim L_i\le 2$ for all $i$, and
satisfying a more stringent but natural narrowness condition
modeled on an analogous condition for pro-$p$ groups.
The two-dimensional homogeneous components of $L$, which include $L_1$, are named {\em diamonds}.
Infinite-dimensional thin Lie algebras with various diamond patterns have been produced,
over fields of positive characteristic, as {\em loop algebras} of suitable finite-dimensional simple Lie algebras,
of classical or of Cartan type depending on the location of the second diamond.
The goal of this paper is a description of the initial structure of a thin Lie algebra, up to the second diamond.

Specifically, if $L_k$ is the {\em second diamond} of $L$, then the quotient $L/L^k$
is a graded Lie algebras of maximal class.
In odd characteristic $p$, the quotient $L/L^k$ is known to be metabelian, and hence uniquely determined up to isomorphism
by its dimension $k$, which ranges in an explicitly known set of possible values:
$3$, $5$, a power of $p$, or one less than twice a power of $p$.
However, the quotient $L/L^k$ need not be metabelian in characteristic two.
We describe here all the possibilities for $L/L^k$ up to isomorphism.
In particular, we prove that $k+1$ equals a power of two.
\end{abstract}
\maketitle

\section{Introduction}\label{sec:introduction}

According to terminology introduced in~\cite{CMNS}, a {\em thin} Lie algebra
is a graded Lie algebra
\[
L=\bigoplus_{i=1}^{\infty}L_{i},
\]
over a field $\F$, such that $\dim L_1=2$
and the \emph{covering property} holds: for every $i\ge 1$ and for every nonzero element $u \in L_{i}$, one has
$[uL_{1}]=L_{i+1}$.
The covering property has the following equivalent formulation (see~\cite{CMNS}):
every graded ideal $I$ of $L$ is located between two consecutive terms $L^{i}$ of the
lower central series of $L$.
The following consequence of the covering property is important:
$L$ is centerless if infinite-dimensional, and the center of $L$ is its highest nonzero homogeneous component otherwise.

The original motivation for such concept came from finite $p$-groups.
Thin $p$-groups were defined in~\cite{Br} according to a certain lattice-theoretic property,
which turned out to be equivalent to the following:
a $p$-group is {\em thin} if it is $2$-generated
and each of its normal subgroups is located between consecutive
terms of the lower central series.
The definition extends naturally to pro-$p$ groups, where a still mysterious example
is the so-called Nottingham group
(see~\cite{Cam} or~\cite{L-GMcKay}, and~\cite{Ershov:finitely_presented} for more recent developments).
In terms of standard invariants
introduced in~\cite{KL-GP},
thin ($p$- or pro-$p$) groups can be described as
groups of {\em obliquity} zero and {\em width} two.

The fact that thinness of a group can be detected, and conveniently exploited as in~\cite{CMNS},
at the level of the graded Lie ring associated with its lower central series
(which is actually a Lie algebra over $\F_p$),
makes it natural to study thin Lie algebras as defined above.
A wealth of examples have been studied and a theory is beginning to emerge,
as one can gauge
from the surveys given in the Introductions of~\cite{CaMa:Hamiltonian} and~\cite{AviMat:A-Z}.

The present paper deals with a piece of the theory which was settled in~\cite{CaJu:quotients}
in a relatively straightforward way for characteristic not two,
but where the case of characteristic two presents serious complications.
In order to motivate and describe the problem we need a broader discussion of thin Lie algebras.

The definition of a thin Lie algebra implies at once that every homogeneous component
of $L$ has dimension one or two (or possibly zero, if $L$ has finite dimension).
A homogeneous component of dimension two is called a {\em diamond} of $L$,
a term reminiscent of the original lattice-theoretic characterization of thin groups in~\cite{Br}.
Thus, $L_1$ is a diamond, and the definition allows that this be the only one.
However, in that case $L$ is a {\em graded Lie algebra of maximal class}, as introduced in~\cite{CMN}.
(Note that the $j$th Lie power $L^j$ of a thin Lie algebra $L$, that is, the $j$th term of its lower central series,
is the homogeneous ideal $\bigoplus_{i\ge j}L_i$, and hence $\dim(L/L^j)=\sum_{i=1}^{j-1}\dim L_i$.)
Because those algebras have been thoroughly studied in~\cite{CN} and~\cite{Ju:maximal},
up to a complete classification of the infinite-dimensional ones (see Section~\ref{sec:main} for more details),
we conveniently exclude them from the definition of a thin Lie algebra.

Thus, a thin Lie algebra $L$ has at least one diamond besides $L_1$,
and the earliest in order of occurrence, say $L_k$, is the {\em second diamond} of $L$.
The structure of the graded Lie algebra of maximal class $L/L^k$ is crucial in the study of $L$.
When the characteristic is not two, and assuming $L$ infinite-dimensional, or at least of dimension large enough,
$L/L^k$ is known from~\cite{CaJu:quotients} to be metabelian, and hence (as one easily sees) determined by $k$ up to isomorphism.
(Beware that the formulation of this result in~\cite{CaJu:quotients} is misleading,
as it wrongly states that every quotient of $L$ which has maximal class is metabelian;
however, $L$ has quotients of maximal class of dimension exceeding by one that of $L/L^k$,
and only one of those is metabelian, namely, $L/[L^2,L^2]$.)
On this property relies the proof in~\cite{AviJur} (which extends more specialized arguments in~\cite{CMNS}) that
$k$ can only take the values $3$, $5$, $p^e$ and $2p^e-1$, where $p$ is the characteristic, if positive.
On the basis of this result one can conveniently divide the thin Lie algebras into three classes:
the {\em classical} ones with $k=3$ or $5$, see~\cite{Mat:thin-groups};
those {\em of Nottingham type}, with $k=p^e$,
see~\cite{Car:Nottingham,Car:Zassenhaus-three,Avi,CaMa:Nottingham,AviMat:A-Z};
those with $k=2p^e-1$, see~\cite{CaMa:thin,CaMa:Hamiltonian,AviMat:-1}.
There may clearly be overlaps between the first case and one of the other two in small characteristics.

Because thin Lie algebras with second diamond in degree $3$ or $5$ are well understood
(apart from one subclass which has only been partially investigated in~\cite{GMY}),
and because only those values can occur in characteristic zero,
the main interest lies now in the modular case, and we will assume from now on
that the ground field $\F$ has positive characteristic $p$.
Also, some thin Lie algebras certainly arise as the graded Lie algebras
(associated with the lower central series)
of thin (pro-) $p$-groups,
as this was the original motivation for their study.
However, those form a small minority among all thin Lie algebras.
In particular, only $3$, $5$ and $p$ occur as possible degrees of the second diamond
for thin algebras associated with groups, the last one being the case for the Nottingham group.
Nevertheless, thin Lie algebras associated with groups have provided initial guidance
in the study of general thin Lie algebras.
For example, the fact that $L/L^k$ must be metabelian,
when $L_k$ is the second diamond of the graded Lie algebra $L$ associated with a thin $p$-group,
reflects a fact discovered by Norman Blackburn
in his pioneering study~\cite{Blackburn} of $p$-groups of maximal class, see~\cite[Hauptsatz~14.6 (a)]{Hup}.
No such guidance can come from groups for thin Lie algebras of characteristic two,
as thin $2$-groups do not exist (or rather, they are groups of maximal class,
if we relax our convention of excluding those from the definition of {\em thin}),
according to~\cite[Satz~11.9 (a)]{Hup}.

In fact, the picture for thin Lie algebras of characteristic two turns out to be
much more complex than for odd characteristic.
Such peculiarity of the characteristic two, hardly unexpected in any theory of Lie algebras,
can already be seen in the context of graded Lie algebras of maximal class,
where the classification proof had to be carried out separately for the odd and even characteristics,
in~\cite{CN} and~\cite{Ju:maximal}, with various additional complications in the latter.
In the broader context of thin Lie algebras, one additional complication of characteristic two
is that certain graded Lie algebras of maximal class may actually be viewed as thin Lie algebras,
the absence of visible diamonds being explained with the presence of {\em fake} diamonds, see Section~\ref{sec:main_thin}
for the meaning of these terms.
Closer to the main object of the present paper, if $L_k$ is the second diamond of a thin Lie algebra in characteristic two
(as we assume without further mention in the rest of this Introduction),
$L/L^k$ need not be metabelian.
A family of infinite-dimensional exceptions was constructed in~\cite{Ju:quotients}, and it was claimed in that paper
that those were the only exceptions.
Unfortunately, the proof given there is incorrect.

In fact, further examples of infinite-dimensional thin Lie algebras
in characteristic two with $L/L^k$ not metabelian, encompassing those given in~\cite{Ju:quotients} as special cases,
were constructed in~\cite{Young:thesis}.
The failure of the statement that $L/L^k$ be metabelian in characteristic two
is so bad that Young's method in~\cite{Young:thesis} produces uncountably many (pairwise not isomorphic) counterexamples
(even over $\F_2$).
This is, roughly speaking, because Young's thin Lie algebras are obtained starting from a certain subclass
of the graded Lie algebras of maximal class, of which uncountably many can be produced
through a recursive procedure and taking suitable limits as in~\cite{CMN}.
As if this were not complex enough, further examples
of thin Lie algebras in characteristic two with $L/L^k$ not metabelian
were built in~\cite{AviMat:A-Z}, as {\em loop algebras}
of certain finite-dimensional simple Lie algebras of Cartan type.

Despite the richness of examples, the structure of the quotient $L/L^k$ for thin Lie algebras in characteristic two
still admits a precise description, which is the main goal of the present paper.
We refer the reader to Theorem~\ref{thm:main} for the details,
but here we mention its remarkable consequence that the degree $k$ of the second diamond
must be one less than a power of two, if $\dim L$ is large enough.
This should be compared with the easier analogue in odd characteristic $p$ mentioned above,
where $k$ can have the form $3$, $5$, $p^e$ or $2p^e-1$ for some $e\ge 1$,
see Corollary~\ref{cor:second_diamond}.

We should mention at this point that a preliminary draft of work in preparation~\cite{JuYo:quotients}
was cited several times in the literature on thin Lie algebras.
Its declared aim was to correct the errors in~\cite{Ju:quotients}
and provide correct information on the structure of $L/L^k$ in characteristic two, as we do here.
However, that tentative draft was never completed, and the more successful approach taken here
makes it obsolete.
We make no attempt to locate the errors in the proofs in~\cite{Ju:quotients}.
Our proofs are independent of that paper, and organized differently.
In fact, our strategy of proof offers several advantages with respect to that in~\cite{Ju:quotients},
besides avoiding the omissions and mistakes made there.
The most notable is that we obtain explicit upper bounds on $\dim L$
in the various steps if $L$ does not satisfied the desired conclusions,
while~\cite{Ju:quotients} had the blanket assumption that $L$ has infinite dimension.
A similar special care when working in finite dimension is required
when dealing with parts of the theory of graded Lie algebras of maximal class,
which we adapt from~\cite{CN,Ju:maximal} to our present needs in Section~\ref{sec:preparatory}.

The paper is organized as follows.
We state our main result on the structure of $L/L^k$ as Theorem~\ref{thm:main},
after recalling some notions and basic results from the theory of graded Lie algebras of maximal class.
In Section~\ref{sec:main_thin} we explain how graded Lie algebras of maximal class can sometimes
be interpreted as thin Lie algebras of Nottingham type with fake second diamond.
We then recast Theorem~\ref{thm:main} in that language, and supplement it with information
on $L$ past the diamond $L_k$ which follows from our proof.
In Section~\ref{sec:preparatory} we collect and adapt to finite dimension
some deeper facts about graded Lie algebras of maximal class
needed in our proof of Theorem~\ref{thm:main}.
Finally, that proof occupies the entire Section~\ref{sec:proof}.
We have placed suitable comments in our proof so that one can also read off
the much easier case of odd characteristic, which was originally dealt with in~\cite{CaJu:quotients}.

Extensive machine calculations performed with
{\sf GAP}~\cite{GAP} have provided invaluable hints for our strategy of proof
in the early stages of this research, and correctness checks in later stages.

\section{Main result}\label{sec:main}

As outlined in the Introduction, our main goal is a description of $L/L^k$
for $L$ a thin Lie algebra with second diamond $L_k$, in the exceptional case where $L/L^k$ is not metabelian,
that is, $[L^2L^2]\not\subseteq L^k$.
According to~\cite{CaJu:quotients}, if $\dim L$ is large enough this can only occur in characteristic two.
Because $L/L^k$ is a graded Lie algebra of maximal class as in~\cite{CMN},
we start by recalling some notions and results from the theory of graded
Lie algebras of maximal class initiated in~\cite{CMN} and further developed in~\cite{CN}.

A {\em graded Lie algebra of maximal class} is a Lie algebra
\[
M=\bigoplus_{i=1}^{\infty}M_{i},
\]
over a field $\F$, graded over the positive integers, with $\dim M_1=2$,
$\dim M_i\le 1$ otherwise, and $[M_i,M_1]=M_{i+1}$.
The notation allows for finite-dimensional algebras, that is, some homogeneous component $M_i$ may be the zero space,
along with all the subsequent homogeneous component.
Note that the $j$th Lie power of $M$ or, in different language, the $j$th term of the lower central series of $M$,
is simply $M^j=\bigoplus_{i\ge j}M_i$.
Hence if $M_{j-1}\neq 0$ then $M/M^j$ has dimension $j$ and nilpotence class $j-1$, hence maximal
with respect to its dimension, thus justifying the name.

In characteristic zero, graded Lie algebras of maximal class are uninteresting, because
the quotient modulo their center turns out to always be metabelian, and hence
uniquely determined up to isomorphism by its dimension.
By contrast, graded Lie algebras of maximal class in prime characteristic can be quite complicated:
it was shown in~\cite{CMN} that there are uncountably many isomorphism classes of them.
Nevertheless, at least the infinite-dimensional ones can be satisfactorily classified, if only in a recursive way.
This was achieved in~\cite{CN} and~\cite{Ju:maximal}, for characteristic odd and two, respectively.

The {\em two-step centralizers} of $M$ are the subspaces of $M_{1}$ defined by
\[
C_{i}=C_{M_{1}}(M_{i})=\{a \in M_{1} : [ab]=0\text{ for all }b \in M_{i}\}.
\]
They owe their name to an analogous concept from Blackburn's theory of $p$-groups of maximal class~\cite{Blackburn},
but for graded Lie algebras of maximal class they are even more crucial,
as their sequence $\{C_{i}\}$ determines $M$ up to isomorphism.
Aside from the trivial case where $M=M_1$ we have $C_1=0$, and this is usually dropped from the sequence.
Because $\dim C_i=2-\dim M_{i+1}$ for $i>1$,
we have $C_i=M_1$ from some point on if $M$ is finite-dimensional.

It is traditional, and convenient, to choose a pair $x,y$ of homogeneous generators for $M$ as follows.
We choose $y \in M_{1}$ such that $C_{2}=\F y$.
If $C_{i}=\F y$ for every $i\geq 2$ then $M$ is the unique metabelian Lie algebra of maximal class
of its dimension.
Otherwise, we choose $x\in M_1$ so that $C_{\q}=\F x$ is the second (distinct) two-step centralizer in order of occurrence.
To make this more concrete, it means that a (highly redundant) presentation for the quotient $M/M^{\bar q+2}$
in the variety of nilpotent Lie algebras of class at most $\bar q+1$,
is given by
\begin{align*}
\langle x,y\colon
[yx^iy]=0
\text{ for $i<\bar q-1$, }
[yx^{\bar q}]=0
\rangle.
\end{align*}
Here we use the left-normed convention for long Lie brackets, as in $[uvw]=[[uv]w]$,
and any exponent in such a long Lie bracket denotes repetition of that entry the specified number of times.

It turns out that $\q=2q$, where the {\em parameter} $q=p^e$ is a power of the characteristic of $\F$,
possibly except for some algebras $M$ of dimension not much larger than $\bar q$.
Precise bounds on $\dim M$ in these exceptional cases can be read off~\cite[Lemma~5.4]{CMN},
but the following formulation given in~\cite{Mat:chain_lengths} is more useful for our purposes:
if there is at least another $C_r$ past $C_{\bar q}$ which is different from $C_2$ and $M_1$
(that is, if $M$ has at least two {\em constituents,} whose definition we recall below), then
$\q$ equals twice some power of the characteristic.

Our choice of denoting by $\q$ what is simply $2q$
is more than a cheap shorthand in complex calculations to follow.
In fact, $\bar q$ is itself a power of the characteristic when $\F$ has characteristic two,
and plays the role of principal parameter in another interpretation of certain graded Lie algebras of maximal class,
to which we devote Section~\ref{sec:main_thin}.

It is easily seen (see~\cite[Lemma~3.3]{CMN})
that at least one of each pair of consecutive homogeneous components of $M$ is centralized by $y$.
(This amounts to the fact that $(\ad y)^2=0$, that is, $y$ is a {\em sandwich element;}
see~\cite{Mat:chain_lengths} for the appropriate extension of this fact to thin Lie algebras.)
Consequently, the sequence of two-step centralizers of $M$ generally consists of consecutive occurrences of $\F y$
interrupted by isolated occurrences of different two-step centralizers.
The distances between these isolated occurrences take the name of {\em constituent lengths}.
More precisely, if
\[
C_{i}\neq \F y, \quad C_{i+1}=\cdots=C_{i+r-1}=\F y, \quad C_{i+r}\not\in\{ \F y, M_1\},
\]
then we say that $\{C_{i+1},\ldots, C_{i+r}\}$ is a {\em constituent} of $M$, of {\em length} $r$.
As in~\cite{CN}, it will be occasionally convenient to say that a nonzero element of $L_i$
is {\em at the beginning} of the constituent or, equivalently, {\em at the end} of the previous constituent.

Our definition of constituent is essentially the same as the updated definition in~\cite{CN},
and hence our constituent lengths are increased by one with respect to~\cite{CMN,CaJu:quotients,Ju:quotients}.
The first constituent $\{C_{2}, \ldots C_{\q}\}$ is conventionally said to be of length $\q$.
However, in our definition of constituent we require the last centralizer $C_{i+r}$ to be different from $M_1$,
which is equivalent to requiring that $M_{i+r+1}\neq 0$.
This detail was not present in~\cite{CN}, which dealt with infinite-dimensional algebras only.
In particular, we do not insist on attaching a constituent length to a possible trailing sequence
of nonzero components centralized by $y$ when $M$ is finite-dimensional.
In fact, while the constituent lengths turn out to be severely restricted,
such final sequence
can be arbitrarily shortened by passing to suitable quotients of $M$,
and so the only general claim one can make about it is that it comprises at most $\q-1$ nonzero components.

With this more restrictive definition of constituent, and under the assumption that $M$ has at least two constituents,
the length of any constituent of $M$ can only be $\q$, or $\q-p^s$ for some $0\leq s \leq e$,
see~\cite[Proposition~5.6]{CMN}, or~\cite{Mat:chain_lengths} for a shorter proof.
Note that when $M$ has only two distinct two-step centralizers (that is, $\F y$ and $\F x$),
as in Theorem~\ref{thm:main} below,
the sequence of constituent lengths actually describes $M$ completely if $\dim M=\infty$,
or if $C_{\dim M-2}\neq C_2$.
Note also that the constituent lengths add up to $\dim M-2$ in those cases.

Now let $L$ be a {\em thin} Lie algebra.
As defined at the beginning of the Introduction, $L$ is a graded Lie algebra
\[
L=\bigoplus_{i=1}^{\infty}L_{i},
\]
such that $\dim L_1=2$, and
$[uL_{1}]=L_{i+1}$ for every $i\ge 1$ and for every nonzero element $u \in L_{i}$
(the {\em covering property}).
Each two-dimensional component of $L$ is a {\em diamond,}
and we are assuming in the definition that $L_1$ is not the only diamond,
to exclude the possibility that $L$ is a graded Lie algebra of maximal class.
Thus, if $L_k$ is the second diamond, $L/L^k$ is a graded Lie algebra of maximal class,
and our primary goal is to describe that up to isomorphism.
The description is very simple in odd characteristic, according to the main result of~\cite{CaJu:quotients},
which we quote in a more precise form.

\begin{theorem}[\cite{CaJu:quotients}]\label{thm:odd}
Let $L$ be a thin Lie algebra, over a field of odd characteristic,
with second diamond $L_k$ and dimension
larger than $4k/3-1$.
Then $L/L^k$ is metabelian, that is, $[L^2L^2]\subseteq L^k$.
\end{theorem}

Because a metabelian graded Lie algebra of maximal class
is determined up to isomorphism by its dimension $k$, and the
possible values for $k$ have been determined in~\cite{AviJur} for $p$ odd
(see our Corollary~\ref{cor:second_diamond} below),
the possible isomorphism types for $L/L^k$ are known in odd characteristic.
In this paper we settle the case of characteristic two,
but add enough comments in our proof to allow one to extract the much easier proof of Theorem~\ref{thm:odd} from it.

The quotient $L/L^k$ of $L$ is not the largest which has maximal class.
In fact, the largest such (graded) quotients of $L$ are obtained by taking
$M=L/(U+L^{k+1})$, a central extension of $L/L^k$, where $U$ is any one-dimensional subspace of $L_k$.
Now $M$ will end with a proper constituent as defined above provided $U\neq [L_{k-1}y]$,
where $C_{L_1}(L_2)=\F y$, and so we assume that.
It will turn out that, under the additional hypothesis that $L/L^6$ is metabelian,
which we discuss at the end of this section, there are exactly two distinct two-step centralizers
$\F y$ and $\F x$ among the $C_{L_1}(L_i)$ for $2\le i<k-1$,
and so this statement extends to $M$ provided we choose $U=[L_{k-1}x]$.
With a harmless abuse, we identify
the generators $x,y$ of $L/L^k$ and $M$ with those of $L$,
and hence corresponding two-step centralizers in them.

\begin{theorem}\label{thm:main}
Let $L$ be a thin Lie algebra, over a field of characteristic two,
with second diamond $L_k$ and dimension
larger than $(4k+1)/3$.
Suppose that the quotient $L/L^{k}$ is not metabelian, but $L/L^6$ is.
Let $\F y=C_{L_{1}}(L_{2})$ and $\F x=C_{L_{1}}(L_{\q})$ be the first two distinct
two-step centralizers, at their first occurrence,
and set
$M=L/([L_{k-1}x]+L^{k+1})$.
Then
\begin{enumerate}
\item the graded Lie algebra of maximal class $M$ has exactly two distinct two-step centralizers;
\item the sequence of constituent lengths of $M$ is
either $\q,\q-2$, or
$\q, \q-1, \q^{\ 2r-3},\q-1$, where $\q$ and $r$ are powers of two;
\item $k+1$ is a power of two.
\end{enumerate}
\end{theorem}

Here $\q^{\ 2r-3}$ stands for $2r-3$ consecutive occurrences of a
constituent of length $\q$, a standard notation as in~\cite{CN}.
The third statement of Theorem~\ref{thm:main} follows at once
from the second, as the constituent lengths in the sequence add up to $2r\bar q-2$.

Note that, because $M$ has only two distinct two-step centralizers,
it is completely determined by its sequence of constituent lengths.
In terms of the sequence of two-step centralizers of $M$
the second assertion of Theorem~\ref{thm:main} can be rephrased as follows:
all two-step centralizers of $M$ are equal to $\F y$,
with the exception of $C_{\bar q}=\F x$ and $C_{i\bar q-1}=\F x$ for
$1<i<(k+1)/\bar q$.
(Note that this range for $i$ is empty when $r=(k+1)/(2\bar q)$ equals $1$.)
We will give yet another equivalent formulation of Theorem~\ref{thm:main}
as Theorem~\ref{thm:main_thin} of the next section.

The sequences of constituent lengths for $M$ obtained in Theorem~\ref{thm:main}
all occur in various infinite-dimensional thin Lie algebras constructed in~\cite{Ju:quotients,Young:thesis,AviMat:A-Z}.
More precisely, the pattern $\q,\q-2$ for the constituent lengths of $M$ occurs in the algebras of~\cite{AviMat:A-Z}
the patterns starting with $\q,\q-1$ occur in~\cite{Ju:quotients},
and all of them occur in~\cite{Young:thesis}.
Furthermore, the algebras produced in~\cite{Young:thesis}
are obtained via an invertible procedure starting from arbitrary graded Lie algebras of maximal class
with only two distinct two-step centralizers.
Because the latter are uncountably many (already over the field of two elements) it follows that there are
uncountably many thin Lie algebras with $L/L^k$ not metabelian and structure of $M$ as described in~Theorem~\ref{thm:main},
for any value of the parameter $r$.

Note that the possibility $\q,\q-2$ for the constituent sequence length of $M$
was omitted in~\cite{Ju:quotients}.
More seriously, the main result of that paper erroneously claimed that
the other possibility, for each value of $r$, would occur for a unique infinite-dimensional thin Lie algebra,
the soluble one constructed there.
The examples mentioned in the previous paragraph show that this is far from the case.

A remarkable consequence of Theorem~\ref{thm:main} is the case $p=2$ of the following assertion.

\begin{cor}\label{cor:second_diamond}
The second diamond of a thin Lie algebra $L$ in characteristic $p$, with $L/L^6$ metabelian and $\dim L$ large enough,
always occurs in odd degree of the form $p^n$ or $2p^n-1$, for some $n>0$.
\end{cor}

In fact, it is known from~\cite{AviJur} (but see~\cite{Mat:chain_lengths} for a shorter proof),
that the second diamond $L_k$ of any thin Lie algebra $L$
occurs in degree of the form $3$, $5$, $q$ or $2q-1$, for some power $q$ of the characteristic,
provided $L$ has dimension large enough, and under the assumption that $L/L^k$ is metabelian.
Granted that $L/L^k$ is known from~\cite{CaJu:quotients} to be always metabelian in odd characteristic,
and that a second diamond in degree $3$ or $5$ would entail that $L/L^6$ is not metabelian,
the case $p>2$ of Corollary~\ref{cor:second_diamond} follows.
The assertion, in Corollary~\ref{cor:second_diamond}, that the second diamond can only occur in odd degree is very elementary,
but serves to exclude the case $p^e$ when $p=2$.
Thus, the conclusion of Corollary~\ref{cor:second_diamond} for $p=2$ is that the second diamond occurs in degree
one less than a power of two, as follows from assertion~(3) of Theorem~\ref{thm:main_thin}.
Incidentally, it is also known that, if we remove the assumption that $L/L^6$ is metabelian,
the second diamond of a thin Lie algebras of characteristic two
can certainly occur in dimension $3$, but not in dimension $5$, as is proved in~\cite{Mat:chain_lengths}.

Our hypothesis that $L/L^6$ is metabelian in the above results needs some justification.
Note that $[L^2L^2]\subseteq L^5$ holds in any two-generated Lie algebra, and we are requiring one step stronger than that.
Generally speaking, the study of thin Lie algebras where $L/L^6$ is not metabelian presents various complications in small
characteristics, due to the scarcity of simple relations of low degree satisfied by a pair of generators.
This occurs, in particular, when the possible degrees $p^e$ and $2p^e-1$ for the second diamond
coincide with $3$ or $5$, as in $p^1=5$, $2p^1-1=5$, $p^1=3$ or $2p^1-1=3$, and hence various classes
of thin Lie algebras with quite different origin mix up.

More specific to our present goals, note that $[L^2L^2]=L_{\bar q+1}$ in the setting of Theorem~\ref{thm:main},
and hence the hypothesis that $L/L^6$ is metabelian is equivalent to requiring that $\q>4$.
The constructions in~\cite{Ju:quotients,Young:thesis,AviMat:A-Z} make sense when $\q=4$,
which is the lowest possible value for $\q=2q$,
and show that all the patterns of constituent lengths for $M$ described in Theorem~\ref{thm:main} occur in this case as well.
However, a proof of Theorem~\ref{thm:main} runs into serious obstacles when $\q=4$.
Our proof covers the subcase where
the sequence of constituent lengths of $M$ begins with $4,3$, that is,
if $L^2/L^6$ is not abelian (and hence has nilpotency class two) but $L^2/L^8$ has nilpotency class three.
Thus, under those assumptions one can still conclude that, provided $\dim L$ is large enough,
$M$ has only two distinct two-step centralizers,
and has constituent length sequence $4,3,4^{\ 2r-3},3$, with $r$ a power of two.
In the remaining case where the constituent length sequence of $M$ starts with $4,2$,
further rather involved calculations, which we do not present here,
show that $M$ has only constituents of length $4$ and $2$, and hence is {\em inflated},
in the language of~\cite{CMN,CN,Ju:maximal}.
However, we know of no counterexample to Theorem~\ref{thm:main} in this case.

\section{Thin Lie algebras of Nottingham type}\label{sec:main_thin}

The fact that $L/L^k$ is always metabelian, for a thin Lie algebra of odd characteristic with second diamond $L_k$,
is a concise way of saying that all one-dimensional homogeneous components between the first and second diamond,
with the necessary exception of the last one, have the same two-step centralizer $\F y=C_{L_1}(L_2)$.
This fact generally extends to most one-dimensional homogeneous components of a thin Lie algebras
not immediately preceding a diamond, with some exceptions.
The exceptions, which occur in both classes of thin Lie algebras with second diamond in degree $q$, and $2q-1$,
for some power $q$ of the characteristic,
are naturally resolved by admitting the presence of {\em fake} diamonds.
Roughly speaking, a fake diamond is a one-dimensional homogeneous component occurring in a degree
where there would be good reasons to expect a diamond,
for example because it fits in a regular pattern
(but deeper reasons are also available).
A peculiarity of characteristic two is that the natural prescription for a fake diamond
for thin Lie algebras {\em of Nottingham type,} that is, with second diamond in degree $q$,
allows for the second diamond itself to be fake.
Of course this contrasts with our detecting the second diamond as the next two-dimensional homogeneous component after $L_1$,
and is the source of many complications of thin Lie algebras of characteristic two,
including the one investigated in this paper.
It is even possible to have a thin Lie algebra of Nottingham type, in characteristic two, where {\em all} diamonds are fake,
and hence the algebra is, in fact, a graded Lie algebra of maximal class.
We refer to~\cite[Section~3]{AviMat:A-Z} for a thorough discussion of this phenomenon,
and only recall here a few basic notions and results which will allow us to recast
Theorem~\ref{thm:main} in the language of thin Lie algebras of Nottingham type.

The peculiarity of characteristic two is better appreciated by introducing the case of odd characteristic first.
Thus, consider a thin Lie algebra
$T=\bigoplus_{i=1}^{\infty}T_{i}$,
in odd characteristic $p$,
with second diamond in degree $q$, a power of $p$.
By definition, $T$ is {\em of Nottingham type}.
If we exclude the rather peculiar case $q=p=3$, we have that $\dim T_3=1$,
and hence there exists a nonzero $y$ in $T_1$ such that $[T_2y]=0$.
Because $p$ is odd, the quotient $T/T^q$ is metabelian according to~\cite{CaJu:quotients}.
This means that $y$ centralizes $T_2,\ldots,T_{q-2}$, and implies that the element
$v=[yx^{q-2}]$ spans $T_{q-1}$, where $x$ is any element of $T_1$ outside $\F y$.
Then $[vx]$ and $[vy]$ span the diamond $T_q$.

Standard calculations then allow one to determine the structure of the quotient
$T/T^{q+2}$, that is, the relations between the generators $[vxx]$, $[vxy]$, $[vxy]$ and $[vyx]$ of $T_{q+1}$.
They can be found in~\cite[Section~2]{AviMat:A-Z},
within a more extensive discussion of thin Lie algebras of Nottingham type,
but we briefly recall them in this paragraph for convenience.
The easiest calculation is $0=[yx^{q-3}[xyy]]=[vyy]$.
The covering property then implies that $T_{q+1}$ is spanned by $[vyx]$
and has, therefore, dimension one.
Consequently, both $[vxx]$ and $[vxy]$ are multiples of $[vyx]$.
In fact, the calculation
\begin{equation*}
0=[yx^{(q-1)/2}[yx^{(q-1)/2}]]=
(-1)^{(q-3)/2}\frac{q-1}{2}[vyx]+(-1)^{(q-1)/2}[vxy]
\end{equation*}
shows that
$[vyx]=-2[vxy]$.
(See Section~\ref{sec:preparatory} for how to perform such calculations rapidly.)
Finally, for any $\beta\in\F$ we have
\[
[v,x+\beta y,x+\beta y]
=[vxx]+\beta\bigl([vxy]+[vyx]\bigr)
=[vxx]-\beta[vxy],
\]
and so we may {\em assume} that $[vxx]=0$ by redefining $x$.

The relations we have just found in degree $q+1$, that is, right after the second diamond,
are a special case of a set of relations, depending on one parameter $\mu\in\F\cup\{\infty\}$,
which generally hold in degree $h+1$ if $T_h$ is any further diamond.
In fact, let $T_{h}$ be any diamond of $T$ past the first.
It is proved in~\cite{Mat:chain_lengths} that no two consecutive homogeneous components
can be diamonds, in a thin Lie algebra with $\dim T_3=1$.
Hence $T_{h-1}$ has dimension one, spanned by $w$, say, and so $[wx]$ and $[wy]$ span $T_h$.
Assuming that $T_{h+1}$ is nonzero, it has dimension one as well, and hence  $[wxx]$, $[wxy]$, $[wxy]$ and $[wyx]$
must satisfy three linear dependence relations.
If they satisfy
\[
[wxx]=0,\qquad [wyy]=0, \qquad (1-\mu)[wxy]=\mu [wyx]
\]
for some $\mu \in \F$, then we say that the diamond $T_h$ is {\em of (finite) type $\mu$.}
In particular, the second diamond $T_q$ has type $-1$.
This terminology extends naturally to include the possibility of an {\em infinite type} $\mu=\infty$,
where we read the third relation as $[wxy]=-[wyx]$.

According to the above definition, types zero and one should not really occur.
In fact, a diamond $T_h$ of type zero should satisfy the relations
$[wxx]=0$ and $[wxy]=0$,
while a diamond of type one should satisfy
$[wyx]=0$ and $[wyy]=0$.
Thus $[wx]$ in the former case, and $[wy]$ in the latter, would be a central element, and hence vanish because
of the covering property.
This contradicts our assumption that $T_h$ has dimension two.
It is convenient, however, to extend our definition of diamond type to include types zero and one,
by relaxing the condition that $\dim T_h=2$.
Thus, a diamond of type $\mu$ will be any homogeneous component $T_h$ such that $T_{h-1}$
is spanned by a single element $w$ which satisfies the relations given earlier.
(We may also assume $T_{h+1}\neq\{0\}$ to make $\mu$ unique.)
When $\mu$ equals zero or one we say that the one-dimensional component $T_h$ is a {\em fake diamond} (of the corresponding type),
and speak of a {\em genuine diamond} in the remaining cases, where $\dim T_h=2$.
The introduction of fake diamonds was first motivated by certain thin Lie algebras of Nottingham type
studied in~\cite{CaMa:Nottingham} (but see also~\cite[Section~7]{AviMat:A-Z}), where the (possibly fake) diamonds occur at regular
intervals and their types follow an arithmetic progression (possibly passing through zero and one).
Fake diamonds turn out to help in the description of more complex types of thin Lie algebras as well,
but we refer to~\cite{AviMat:A-Z} for a survey.

We warn the reader about a possible source of confusion.
A diamond $T_h$ of type one occurs when $T_{h-1}=\F w$, $[wy]=0$ (whence $T_h=\F[wx]$) and $[wxx]=0$.
But then one can prove that $[wxyy]=0$, and so $T_{h+1}=\F[wxy]$ can be interpreted as a diamond of type zero.
This ambiguity of interpretation is easily resolved when the diamonds occur at regular distances.
As a rule, they are spaced apart by $q-1$ degrees, as in~\cite{Car:Nottingham,Car:Zassenhaus-three,CaMa:Nottingham,AviMat:A-Z}.
However, some of the algebras in~\cite{Young:thesis} have sequences of diamonds of type one, say,
occurring at regular intervals of $q$ degrees apart.
Generally speaking, the higher distances there can be justified by reinterpreting a diamond of type one
in a given degree as a diamond of type zero in one degree higher.
A similar situation occurs in Theorem~\ref{thm:main_thin} below, and we give preference to type one over zero in our descriptions.

Various constructions in~\cite{Car:Nottingham,Young:thesis,CaMa:Nottingham,AviMat:A-Z}
for thin Lie algebras of Nottingham type of odd characteristic naturally make sense in characteristic two as well,
and usually produce thin Lie algebras once additional central elements are factored out.
However, the second diamond $T_q$, which always has type $-1$, becomes fake in characteristic two,
and hence may actually not be recognized as the second diamond.
In the extreme case of the algebras of~\cite{Car:Nottingham}, which directly generalize the one associated
with the lower central series of the Nottingham group, all diamonds have type $-1$ and hence
become fake in characteristic two, yielding graded Lie algebras of maximal class.
In other cases genuine diamonds remain, but it is clear that the definition of thin Lie algebras of Nottingham type
needs special treatment in characteristic two.

An appropriate way to resolve this was suggested in~\cite[Section~3]{AviMat:A-Z},
and involves redefining the {\em degree of the second diamond,} in a thin Lie algebra $T$ with $\dim T_3=1$, as
one less than the dimension of the largest metabelian quotient of $T$.
This is consistent with the natural definition in odd characteristic,
but in characteristic two it allows for the possibility of a fake diamond of type one in degree $q$ (a power of two).
We will not fully adopt that definition here, as it would conflict with the terminology used so far,
but point out that this ambiguity in what should be considered the second diamond in characteristic two
may be the true source of the phenomenon investigated in this paper.
However, in order to take advantage of the convenient terminology of diamond types,
we will say that a thin Lie algebra $T$ in characteristic two is {\em of Nottingham type with fake second diamond in degree $q$}
if $T/T^{q+1}$ is metabelian and $T/T^{q+2}$ is a graded Lie algebra of maximal class which is not metabelian.
Because $L/L^5$ is metabelian for every two-generated Lie algebra $L$,
these assumptions imply that $q\ge 4$.
Facts from the theory of graded Lie algebras of maximal class
recalled in Section~\ref{sec:main} imply that $q$ is a power of two,
and $T_q$ may be viewed as a fake diamond of type one.
With this terminology at hand we recast Theorem~\ref{thm:main} as follows.

\begin{theorem}\label{thm:main_thin}
Let $L$ be a thin Lie algebra, over a field of characteristic two,
with second (genuine) diamond $L_k$ and dimension larger than $(4k+1)/3$.
Suppose that the quotient $L/L^{k}$ is not metabelian, but $L/L^6$ is.
Let $\F y=C_{L_{1}}(L_{2})$ and $\F x=C_{L_{1}}(L_{\q})$ be the first two distinct
two-step centralizers, at their first occurrence.
View $T$ as of Nottingham type with fake second diamond in degree $\q$.
Then
\begin{enumerate}
\item each homogeneous component $L_i$, for $1<i<k-1$, is centralized by $x$ or $y$;
\item $k=2r\q-1$, where $r \geq 1$ is a power of two;
\item $L$ has fake diamonds of type one in each degree $m \q-1$,
for $1<m<2r$.
\end{enumerate}
\end{theorem}

Note that $\q=2q$ has taken the place of the usual parameter $q$ for thin Lie algebras of Nottingham type.
It is understood in assertion~(3), and in similar assertions later,
that there are no diamonds of type one other than those mentioned, in the range considered.
We have found convenient to consider only one type of fake diamonds in the statement and the proof, namely, those of type one,
but of course an equivalent formulation of assertion~(3) is that $L$ has fake diamonds of type zero in each degree $m \q$,
for $1<m<2r$.

We have already noted in the previous section that the various isomorphism types for the
quotient $L/L^k$ described in Theorem~\ref{thm:main_thin},
all occur in various infinite-dimensional thin Lie algebras constructed
in~\cite{Ju:quotients,Young:thesis,AviMat:A-Z}.
We add here that the second genuine diamond $L_k$ can have any type $\mu\in\F\setminus\{0,1\}$ or $\mu=\infty$
in the algebras with $r=1$ described
in~\cite[Theorem~6.1]{AviMat:A-Z} and~\cite[Theorem~5.1]{AviMat:A-Z}, respectively.
However, all genuine diamonds, including the second diamond $L_k$, have type $\mu=\infty$
in the algebras with $r>1$ of~\cite{Young:thesis}, of which those exhibited in~\cite{Ju:quotients} are a special case.

Our proof actually describes the structure of $L$ a little after the diamond $L_k$ as well, and allows us to supplement
Theorem~\ref{thm:main_thin} with the following assertions.
The examples from~\cite{Ju:quotients,Young:thesis,AviMat:A-Z} mentioned in the previous paragraph
show that all the parameters involved in these assertions are best possible.

\begin{theorem}\label{thm:main_thin_more}
Assume the hypotheses of Theorem~\ref{thm:main_thin}.
Then
\begin{enumerate}
\setcounter{enumi}{3}
\item each homogeneous component $L_i$, for $k<i<3(k-1)/2$, is centralized by $x$ or $y$;
\item $L$ has fake diamonds of type one in each degree $m \q-2$,
for $2r<m<3r$, as long as $m \q-2\le\dim L+3$;
\item if $r>1$ then the second genuine diamond $L_k$ has type $\infty$.
\end{enumerate}
\end{theorem}

\section{Preliminary work on graded Lie algebras of maximal class}\label{sec:preparatory}

In this section we prove some results on graded Lie algebras of maximal class in view
of their application to a suitable quotient of a thin Lie algebra $L$.
For the most part these results are already known from the papers~\cite{CN,Ju:maximal}
on the classification of graded Lie algebras of maximal class,
except that we need to revise them to cover the case of present interest where the algebra has finite dimension.

The following generalization of the Jacobi identity in a Lie algebra
will be used repeatedly and without specific mention
\[
[v[yx^{n}]]=\sum_{i=0}^{n}(-1)^{n} \binom{n}{i}[vx^{i}yx^{n-i}].
\]
Because we work in characteristic two (with the exception of~Lemma~\ref{lemma:three_centralizers})
we will omit the alternating signs $(-1)^{n}$ in such summations.
However, to improve clarity we will, at least initially, keep some minus signs in our calculations
when they arise from application of the Jacobi identity in the form
$[x[yz]]=[xyz]-[xzy]$,
and then drop them as soon as they cause clutter.
The binomial coefficients in the formula are conveniently evaluated modulo $p$ by means of Lucas' Theorem:
\[
\binom{a}{b}\equiv \prod_{i=0}^{m} \binom{a_{i}}{b_{i}} \pmod{p}
\]
where $a=\sum_{i=0}^{m}a_{i}p^{i}$ and $b=\sum_{i=0}^{m}b_{i}p^{i}$ are the $p$-adic expansions
of the nonnegative integers $a$ and $b$.
For example, two recurring instances of the formula are $[v[yx^{n}]]=[vyx^{n}]-[vx^{n}y]$
when $n$ is a power of the characteristic $p$, and
$[v[yx^{n}]]=\sum_{i=0}^{n}[vx^{i}yx^{n-i}]$
when $n+1$ is a power of $p$.

For later use we recall, and adjust to the finite-dimensional case,
a known result on graded Lie algebras of maximal class with more than
two distinct two-step centralizers, which is valid in arbitrary characteristic.

\begin{lemma}\label{lemma:three_centralizers}
Let $M$ be a graded Lie algebra of maximal class with at least two constituents.
With the standard notation introduced in Section~\ref{sec:main},
let $C_2=\F y$ and $C_{\q}=\F x$ be the first and second two-step centralizers at their earliest occurrence, where $\q=2q$.

\begin{enumerate}
\item
Any two-step centralizer other than $\F y$ and $\F x$, which is not at the end of the last constituent of $M$,
is preceded by a constituent of length $q$ and is followed by a constituent of length $q$.
\item
If $M$ has at least three distinct two-step centralizers, all occurring before the last constituent,
then all constituents of $M$ have length $\q$ or $q$.
\end{enumerate}
\end{lemma}

\begin{proof}
Both assertions were proved in~\cite{CN} under the additional hypotheses that $p>2$ and $M$ has infinite dimension.
The former hypothesis was immaterial in the proof, and we have adjusted our statements in order
to dispose of the latter hypothesis.

More specifically, a proof of assertion~(1) was first given in~\cite{CaJu:quotients},
again under the additional hypotheses mentioned above.
Another proof was then given in~\cite[Lemma~3.12]{CN},
based on the {\em specialization of two-step centralizers} introduced in~\cite[Proposition~4.3]{CN}.
That technique was then essential in the proof of assertion~(2) given in~\cite[Step~11]{CN},
together with assertion~(1) and several other ingredients from~\cite{CN}, or the corresponding ones in~\cite{Ju:maximal}
for the case of characteristic two.
One can verify that all those arguments, originally formulated under the blanket assumption that $M$ has infinite dimension,
remain valid here.
\end{proof}

Note that Lemma~\ref{lemma:three_centralizers} contains no assumption on $\dim L$ beyond what is implicit
in the constituents mentioned, with our updated meaning of constituent.
This also applies to Proposition~\ref{prop:chains} below.

Now consider a graded Lie algebra of maximal class $M$, with exactly
two distinct two-step centralizers and at least two constituents.
As recalled in Section~\ref{sec:main}, letting $\F y=C_{2}$ be the first two-step centralizer,
and $\F x=C_{\q}$ the second at its earliest occurrence, we have $\q=2q$ and $q=p^e$ for some $e\geq 1$.
In general, the lengths of the first two constituents exert great influence on the structure of $M$.
According to~\cite[Theorem~5.5]{CMN}, in odd characteristic the length of the second constituent can only be $q$.
In fact, if $\dim M$ is large enough the first constituent is followed by at least $p-2$ constituents of length $q$.
However, in characteristic two the second constituent can take any length permitted by the general result
on constituent lengths recalled in Section~2, hence $\q$ or $\q-2^s$ for some $0\le s\le e$,
except for the highest value $\q$, which is easily excluded by computing
\[
0=[yx^{\q-1}[yx^{\q-1}]]=\sum_{i=0}^{\q-1}[yx^{\q-1+i}yx^{\q-1-i}]=[yx^{\q-1}yx^{\q-1}].
\]

We will be interested in a situation where the second constituent of $M$, in characteristic two,
attains its highest possible length, $\q-1$.
Then it turns out that any other constituent of $M$ has length $\q$ or $\q-1$, if $M$ has infinite dimension.
This is a rather special case of Step~6 in Jurman's classification~\cite{Ju:maximal}
(anticipated in the Proposition in~\cite[Section~6]{Ju:quotients}),
which is analogous to the corresponding step in~\cite{CN}.
However, we give here a revised proof of this special case which includes the possibility,
neglected in~\cite{Ju:maximal}, that $M$ has finite dimension.

\begin{prop}\label{prop:chains}
Let $M$ be a graded Lie algebra of maximal class
with first constituent of length $\q$, and second constituent of length $\q-1$ (hence in characteristic two).
Then every constituent of $M$ has length $\q$ or $\q-1$ except, possibly, for the last.
However, if the last constituent is shorter than $\q-1$ then the penultimate constituent has length $\q-1$.
\end{prop}

\begin{proof}
According to Lemma~\ref{lemma:three_centralizers},
$M$ has only two distinct two-step centralizers, $\F y$ and $\F x$.
Consequently, of all the iterated Lie brackets $[yxz_3\cdots z_r]$ with $z_i\in\{x,y\}$, for a given $r<\dim M$,
exactly one is nonzero.
We may view this as the {\em standard} choice of a nonzero element of $M_{r}$.
More importantly, keeping this in mind allows one to see at a glance that most such
expressions vanish in the calculations to follow.
This crucial fact can also be expressed by saying that $M$ is actually $\Z^2$-graded,
by assigning independent degrees to $x$ and $y$.

The results of \cite{CMN} imply that each constituent of $M$ has length $h$ of the form $\q$ or $\q-2^s$,
for $0\leq s\leq e$.
Note that $\q-1$ is the only odd number among those.
Hence it suffices to prove that no constituent of $M$ can have even length less than $\q$, with the stated exception.
Set $v_{1}=[yx^{\q-2}]$ and $v_{2}=[v_{1}xyx^{\q-3}]$, so that $v_{1}$ and $v_{2}$ are nonzero elements of degree
$\q-1$ and $2\q-2$ respectively.
As a preliminary observation, note that
$[z[v_1x]]=\sum_{i=0}^{\q-1}[zx^iyx^{\q-1-i}]$ for any $z\in L$.
Also, if $z$ is the standard homogeneous element of $L$ of degree $r$,
then at most one term of the sum is nonzero.
That happens when exactly one of $C_r,\ldots,C_{r+q-1}$ equals $\F x$, and in that case $[z[v_1x]]$ equals
the standard element of degree $r+q$.

Suppose for a contradiction that $M$ has a third constituent, of even length $h<\q$.
Setting $\ell=(h-2)/2$, we obtain the desired contradiction by computing
\begin{align*}
0&=[v_{1}xyx^{\ell}[v_{1}xyx^{\ell}]]
\\&=
\sum_{i=0}^{\ell}\binom{\ell}{i}
[v_{1}xyx^{\ell+i}[v_{1}xy]x^{\ell-i}]
\\&=
[v_{1}xyx^{h-2}[v_{1}xy]]
\\&=
[v_{1}xyx^{h-2}[v_1x]y]-[v_{1}xyx^{h-2}y[v_1x]]
\\&=
[v_{2}xyx^{h-1}y].
\end{align*}
In this calculation, all terms with $i<\ell$ in the summation vanish, because of our assumption that $x$ centralizes $M_{2\q-1+h}$.
We will use similar arguments below without further mention.

Now work inductively and consider a constituent of $M$ after the second one.
Thus, we have $C_{r}=\F x$,
$C_{r+1}=\cdots =C_{r+m-1}=\F y$ and $C_{r+m}=\F x$,
for some $r\geq 3\q -2$,
and we may assume that $m$ equals either $\q-1$ or $\q$.
We will prove that the next constituent has also length at least $\q-1$, except for the case detailed in the statement.
Hence suppose that it has even length $h<\q$.
Let $w$ be a nonzero element in $M_{r-1}$, say the standard element.

Suppose first that $m=\q$.
We have
\begin{equation}\label{eq:chain}
\begin{split}
0&=[wx[v_{1}xyx^{h-1}y]]
\\&=
[wx[v_{1}xyx^{h-1}]y]-[wxy[v_{1}xyx^{h-1}]]
\\&=
[wx[v_{1}xy]x^{h-1}y]-\sum_{i=0}^{h-1}\binom{h-1}{i}[wxyx^{i}[v_{1}xy]x^{h-1-i}]
\\&=
[wx[v_{1}xy]x^{h-1}y]-[wxyx^{h-1}[v_{1}xy]]
\\&=
[wx[v_{1}x]yx^{h-1}y]-[wxy[v_{1}x]x^{h-1}y]
-[wxyx^{h-1}[v_{1}x]y].
\end{split}
\end{equation}
Each of the three summands in the last expression equals $[wxyx^{\q-1}yx^{h-1}y]$,
which is the standard element of its degree,
and we obtain the desired contradiction.

Now let $m=\q-1$.
Note in passing that the above calculation remains valid provided we replace $h-1$ with $h-2$ throughout, and hence
\[
0=[wx[v_{1}xyx^{h-2}y]]
=
[wx[v_{1}x]yx^{h-2}y]-[wxy[v_{1}x]x^{h-2}y]
-[wxyx^{h-2}[v_{1}x]y].
\]
But here the first summand vanishes, while each of the other two summands equals
$[wxyx^{\q-2}yx^{h-1}y]$, and so we cannot obtain a contradiction in this way
(see Remark~\ref{rem:BiZas} below).
However, we can assume that the constituent of length $h$ under examination is not the last constituent of $M$,
and then the calculation
\begin{align*}
0&=[wx[v_{1}xyx^{h-1}y]]
\\&=
[wx[v_{1}xyx^{h-1}]y]-[wxy[v_{1}xyx^{h-1}]]
\\&=
\sum_{i=0}^{h}\binom{h-1}{i}[wxyx^{i}[v_{1}xy]x^{h-1-i}]
\\&=
(h-1)[wxyx^{h-2}[v_{1}xy]x]
\\&=
[wxyx^{h-2}[v_{1}x]yx]-[wxyx^{h-2}y[v_{1}x]x]
\\&=
[wxyx^{\q-2}yx^{h-1}yx]
\end{align*}
provides the desired contradiction.
\end{proof}

\begin{rem}\label{rem:BiZas}
In the exceptional case of Proposition~\ref{prop:chains}, if the penultimate constituent of $M$ has length $\q-1$
then the last constituent can indeed have any length
allowed by the general theory of constituent lengths.
In fact, the lengths $\q-2^s$ with $0<s\leq e$
occur by taking central extensions of suitable quotients of the infinite-dimensional graded Lie algebras
of maximal class $B_\ell(g,h)$ produced in~\cite{Ju:Bi-Zassenhaus}.
The algebras $B_\ell(g,h)$ are obtained as {\em loop algebras} of certain finite-dimensional simple Lie algebras $B(g,h)$
(where $g,h$ are integral parameters with $g\ge 2$ and $h\ge 1$) with respect to a suitable nonsingular derivation.
The algebra $B_\ell(g,h)$ has periodic constituent length sequence $\q,(\q-1,\q^{2^{g}-2},\q-1)^\infty$, where $\q=2^{h+1}$
(and so $e=h$ here).
Hence any (graded) quotient of $B_\ell(g,h)$ has only constituents $\q$ and $\q-1$.
In particular, the quotient of dimension $2\q+2^g\q-1$
has sequence of constituent lengths
$\q,\q-1,\q^{2^{g}-2},\q-1,\q-1$.
However, because the underlying simple Lie algebra $B(g,h)$ has non-trivial central extensions,
certain quotients of $B_\ell(g,h)$ have extensions by a one-dimensional center, hence themselves of maximal class,
with sequence of constituent lengths of the form $\q,\q-1,\q^{2^{g}-2},\q-1,\q-2^s$, for any $s$ with with $0<s\leq h$.
Specifically, such algebra is obtained as a central extension of the quotient of $B_\ell(g,h)$
of dimension $2\q+2^g\q-2^s-1$ by a central element arising from the cocycle of $B(g,h)$ denoted by
$\psi_{s+g-1}$ in~\cite[Section~4]{Ju:Bi-Zassenhaus}.
In fact, the same cocycles are responsible for central extensions of the quotients of $B_\ell(g,h)$
of dimension $2\q+1-2^s$, which have sequence of constituent lengths $\q,\q-2^s$.
\end{rem}

Note that, for constituents later than the third, the proof of Proposition~\ref{prop:chains}
uses only the initial structure of $M$, say about its first two constituents,
and hence its inductive steps may work also for algebras which may not be of maximal class,
but which have a suitable quotient of maximal class.
In particular, calculations similar to Equation~\eqref{eq:chain} will play a role in establishing that fake
diamonds of type one occur at regular distances of $\q$ in the thin Lie algebra $L$
of Subsection~\ref{subsec:t=2q-2}.

\section{Proof of Theorem~\ref{thm:main}}\label{sec:proof}

Let $L$ be a thin Lie algebra with second diamond $L_k$ and $L/L^k$ not metabelian.
It was shown in~\cite{CaJu:quotients} how in odd characteristic it follows that $L$ has finite dimension.
The proof in~\cite{CaJu:quotients} was based on considering the slightly larger quotient
$M=L/([L_{k-1}x]+L^{k+1})$ of maximal class of $L$,
where $\F x$ is the second two-step centralizer of $L/L^k$, as usual.
A contradiction or a bound on $\dim L$ was then obtained,
according to each possible length of the last constituent of $M$, denoted there by $t+1$.
Specifically, that proof showed that $t+1=2q$ cannot occur at all, that $\dim L\le k+2$ if $t+1<2q-1$,
and that $\dim L\le k+q+1$ if $t+1=2q-1$.
We proceed similarly here in characteristic two, except that here no bound on $\dim L$ exists when $t+1$ equals $\q-1$ or $\q-2$.
With considerably more effort we will be able to determine the structure of $M$ in those cases,
provided $\dim L$ is large enough.

We divide the proof into several parts for clarity, which we present in corresponding subsections.
Assume that $\dim L$ is large enough in the following brief description.
In Subsection~\ref{subsec:initial} we study the structure of $L$ a little before and a little after the second diamond.
In particular, we show that the last constituent of the quotient $M=L/([L_{k-1}x]+L^{k+1})$ of $L$
can only have length $\q-2$ or $\q-1$.
In Subsection~\ref{subsec:second} we prove that if the second constituent of $M$ is shorter than $\q-1$,
which is the highest length allowed by the general theory,
then $M$ has constituent length sequence $\q,\q-2$.
This is one of the possibilities allowed by Theorem~\ref{thm:main}, and so we may now assume that
the second constituent of $M$ has length $\q-1$.
Then we show in Subsection~\ref{subsec:t=2q-3} that the last constituent of $M$ must have length $\q-1$ as well.
Proposition~\ref{prop:chains} now implies that the constituents of $M$ can only have length $\q-1$ or $\q$.
In Subsection~\ref{subsec:final}, which is the most complex part of the proof,
we show that a further constituent of $M$ of length $\q-1$ besides the second and the last
leads to a contradiction.
Hence the sequence of constituent lengths of $M$ has the form $\q, \q-1, \q^n,\q-1$ for some $n$.
Finally, in Subection~\ref{subsec:t=2q-2} we prove that $n+3$ is a power of two.
We comment in Subsection~\ref{subsec:t=2q-2} on the reasons for reversing the order of these last two parts of the proof.

\subsection{Initial considerations}\label{subsec:initial}
We can naturally identify
$M_r$ with $L_r$ for $r<k$.
In particular, the spaces $C_r=C_{L_1}(L_r)$ for $r<k-1$ coincide with the two-step centralizers  $C_r=C_{M_1}(M_r)$ of $M$,
but we complete their list with $C_{k-1}=\F x=C_{M_1}(M_{k-1})$, rather than with $C_{L_1}(L_{k-1})=\{0\}$.
Thus, the spaces $C_r$ are the two-step centralizers of the graded Lie algebra of maximal class $M$.
Because we are assuming that $L/L^k$ is not metabelian, $M$ has at least two constituents, and hence
the length $\q$ of the first one must be (twice) a power of two, $\q=2\cdot 2^e>2$.
Furthermore, the possible lengths of any constituent of $M$ are $\q$ and $\q-2^s$ for some $0\leq s\leq e$.

We can say more about the length of the last constituent of $M$ if we look at $L$ rather than just its quotient $M$,
assuming that $L_{k+1}\neq\{0\}$.
The following result relies on the same arguments used in~\cite{CaJu:quotients} for the case of odd characteristic,
which we quote for completeness.
Before that we fix some notation which will be kept through the whole section.
Let $t+1$ be the length of the last constituent of $M$, and
choose a nonzero element $u$ of $L_{k-t-3}$.
Then $[uy]=0$ and $L_{k-t-2}=\F [ux]$, but $C_{k-t-2}=\F(x-cy)$ for some $c\in\F$, and so $[uxx]=c[uxy]$.
Furthermore, $[uxyx^{i}y]=0$ for $0\le i\le t-1$
and hence $L_{k-t-1+i}=\F [uxyx^{i}]$ for $0\le i\le t$.
We also set $v=[uxyx^{t}]$, so that $[vx]$ and $[vy]$ span the diamond $L_k$.

\begin{lemma}\label{lemma:1}
If $\dim L\ge k+3$, the last constituent of $M$ has length $\q-2$ or $\q-1$.
\end{lemma}

\begin{proof}
The possibility that $t+1=\q$ is easily ruled out by computing
\[
0=[ux[yx^{\q}]]=[uxyx^{\q}]-[ux x^{\q}y]=[vx]-c[vy].
\]
We can also deal with the cases $t+1<\q-2$ as done in~\cite{CaJu:quotients} in odd characteristic.
Thus, after noting that $0=[uxyx^{t-1}[xyy]]=[vyy]$, whence $L_{k+1}=\F [vyx]$ because of the covering property, we compute
\begin{align*}
0&=[u[yx^{t+2}y]]=[u[yx^{t+2}]y]
\\&=
\sum_{i=0}^{t+2}\binom{t+2}{i}[ux^{i}yx^{t+2-i}y]
\\&=
(t+2)[uxyx^{t+1}y]+[ux^{t+2}yy]
\\&=
(t+2)[vxy]+c[vyy]=[vxy],
\end{align*}
as $t$ is odd.
This is a contradiction.
\end{proof}

Note that the analogous steps in odd characteristic would leave only one possibility for $t+1$, namely, $\q-1$;
this possibility was then excluded in~\cite{CaJu:quotients},
in a way which we recall for completeness after the proof of our Lemma~\ref{lemma:2}.
In characteristic two, neither of the values $\q-1$ and $\q-2$ for $t+1$ implies that $L$ has finite dimension.

From now on we assume without mention that $\dim L\ge k+3$,
which is just one step further than the implicit assumption that $L_k$ is a genuine diamond,
and is precisely enough to assign $L_k$ a diamond type.
Stronger assumptions on $\dim L$ will be specified as we need them in the individual arguments,
and will always be implied by the hypothesis $\dim L>(4k+1)/3$ of Theorem~\ref{thm:main}.

We can strengthen the conclusion of Lemma~\ref{lemma:1} in the special case where $M$ has only two constituents.
In fact, the same calculation
$0=[yx^{\q-1}[yx^{\q-1}]]=[vx]$
used in section~\ref{sec:preparatory} to to establish the upper bound $\q-1$ for the length of the second constituent
of an arbitrary graded Lie algebra of maximal class,
shows in the present situation that the second constituent of $M$ can only have length $\q-2$.

Now we are ready to obtain some information on $L$ past the diamond $L_k$, but we need the crucial assumption
that $M$ has only two distinct two-step centralizers.
According to Lemma~\ref{lemma:three_centralizers}, this holds under our hypothesis that $\q>4$,
because then $t+1$, which Lemma~\ref{lemma:1} allows to be only $\q-1$ or $\q-2$, is different from both $\q=2q$ and $q$.
In fact, Lemma~\ref{lemma:three_centralizers} guarantees that $M$ has only two distinct two-step centralizers
also when $\q=4$, provided $M$ has at least one constituent of length $3$ (hence different from $\q$ and $q$).
In particular, if the sequence of constituent lengths of $M$ begins with $4,3$, then
our arguments in Subsections~\ref{subsec:t=2q-3}, \ref{subsec:t=2q-2} and~\ref{subsec:final}
will remain valid and prove the conclusion of Theorem~\ref{thm:main} also in this case,
as mentioned at the end of Section~\ref{sec:main}.

\begin{lemma}\label{lemma:2}
Suppose that $M$ has only two distinct two-step centralizers.
Then $[L_{k-1}xx]=0$.
Furthermore, $L_{k+1},\ldots,L_{k+q-1}$ are centralized by $y$.
\end{lemma}

\begin{proof}
With notation as above we have $L_{k-1}=\F v$.
We have already noted in the proof of Lemma~\ref{lemma:1} that $[vyy]=0$, and hence $[vyx]$ is nonzero and spans $L_{k+1}$.
We now prove that $[vxx]=0$, using the fact that $[uxx]=0$ because $[uxy]$ is nonzero and $M$ has only
two distinct two-step centralizers.
This is slightly easier when $t+1=\q-1$, as
\[
0=[ux[yx^{\q}]]=[uxyx^{\q}]=[vxx].
\]
When $t+1=\q-2$ we only obtain that
\[
0=[ux[yx^{\q}]]=[uxyx^{\q}]=[vxxx].
\]
Suppose for a contradiction that $[vxx]\neq 0$, whence
$[vxx]$ spans $L_{k+1}$, and the latter is centralized by $x$.
In particular, $[vxyx]=0$,
but also
\[
0=[uxy[yx^{\q-2}y]]=[uxy[yx^{\q-2}]y]=[vxyy].
\]
Therefore, $[vxy]$ is central, and hence zero by the covering property (provided $\dim L\ge k+4$).
If $[vxx]=a[vyx]$ then $[vx]-a[vy]$ is also central, and hence $[vx]=a[vy]$.
This contradicts the fact that $L_{k}$ has dimension two,
and we conclude that $[vxx]=0$ as desired.

An argument from~\cite{CaJu:quotients} applies unchanged in characteristic two to prove
our second assertion, which is analogous to the fact that $q$ is a lower bound for the length of a constituent
in graded Lie algebras of maximal class.
In fact, because $t+1\geq q$ we may consider the element $[uxyx^{t+1-q}]$.
Recall that $[uxyx^{t}]=v$
and note that $[uxyx^{t+1-q+j}y]=0$ for $0\le j<q-1$.
Because $[yx^{q+j}y]=0$ in the same range we have
\begin{align*}
0&=[uxyx^{t+1-q}[yx^{q+j}y]]
=[uxyx^{t+1-q}[yx^{q+j}]y]
\\&=
\binom{q+j}{q-1}[uxyx^{t}yx^{j+1}y]+\binom{q+j}{q}[uxyx^{t+1}yx^{j}y]
=[vxyx^{j}y],
\end{align*}
where we have used the fact that $[vxx]=0$.
Because of the covering property it follows inductively that
$L_{k+i+1}=\F [vxyx^{i}]$ for $0\le i<q$, and the desired conclusion follows.
\end{proof}

For comparison, we quote from~\cite{CaJu:quotients} the conclusion of the proof in the case of odd characteristic,
after noting that the analogue of Lemma~\ref{lemma:2} holds there as well, and with a simpler proof.
As pointed out earlier, the analogue of Lemma~\ref{lemma:1} leaves us with $t+1=\q-1$,
whence $v=[uxyx^{\q-2}]$.
Then the calculation
\[
[uxyx^{q-2}[yx^{\q}]]
=\binom{\q}{q+1}[vxyx^{q-1}]
-\binom{\q}{q}[vyx^q]
=-2[vyx^q]
\]
together with Lemma~\ref{lemma:2} implies that $\dim L\le k+q+1$.
Theorem~\ref{thm:odd} easily follows from here, but
the coefficient $-2$ makes this argument inconclusive in characteristic two.

According to Lemma~\ref{lemma:2}, if $M$ has only two distinct two step centralizers, and $\dim L\ge k+q+2$,
then the diamond $L_{k}$ is followed by at least
$q$ one-dimensional homogeneous components,
all centralized by $y$ with the possible exception of the last.
The fact that $[L_{k-1}xx]=0$ implies that each of $[vxy]$ and $[vyx]$ spans $L_{k+1}$, because of the covering property.
Furthermore, the diamond $L_k$ can be assigned a type $\mu$
as in Section~\ref{sec:main_thin}, which is determined by
$[wyx]=(\mu^{-1}-1)[wxy]$.
As noted after Theorem~\ref{thm:main}, the infinite-dimensional thin Lie algebras
constructed in~\cite{AviMat:A-Z} show that, in general,
the coefficient $\mu^{-1}-1$ may take any value in $\F\setminus\{0\}$.

\subsection{Second constituent of $M$ shorter than $\q-1$}\label{subsec:second}

Suppose that the second constituent of $M$ is shorter than it highest possible value $\q-1$.
Our goal in this case is to prove that the quotient $M$ of $L$ has only two constituents,
and that $M$ has constituent length sequence $\q,\q-2$.
As noted after Theorem~\ref{thm:main}, this configuration occurs for several thin Lie algebras
constructed in~\cite{Young:thesis,AviMat:A-Z}.

According to Lemma~\ref{lemma:1}, the last constituent of $M$ has length $\q-2$ or $\q-1$.
If $\dim L$ is large enough, we prove the desired conclusion in the former case,
and we show that the latter case leads to a contradiction.

In the proofs we need to consider not only the last constituent of $M$, but also the previous constituent.
Thus, let $r+1$ be the length of the penultimate constituent of $M$.
If $M$ has at least three constituents, the constituent length sequence of $M$
ends with $r+1,t+1$, and $L_{k-t-r-4}$ is one-dimensional, spanned by an element $w$.
Because $M$ has only two distinct two-step centralizers,
the relations which hold in $L$ between $\F w$ and the diamond $L_k$ can be conveniently formulated as follows.
All elements of the form $[wz_1\ldots z_{t+r+3}]$ with $z_i\in\{x,y\}$ vanish, with the only
exception of $[wxyx^ryx^t]$, which equals $v$ after replacing $w$ with a scalar multiple as we may.

\begin{lemma}\label{lemma:3}
Assume $\q>4$.
Suppose that the second constituent of $M$ is shorter than $\q-1$, and that the last constituent of $M$ has length $\q-2$.
Then $M$ has constituent length sequence $\q,\q-2$, unless $\dim L\le k+q+2$.
\end{lemma}

\begin{proof}
We may assume that $M$ has more than two constituents, and hence that as$[v_1xyx^{\q-2}]$ vanishes in $L$
(and not just in $M$), where we have conveniently set $v_1=[yx^{\q-2}]$ as in the proof of Proposition~\ref{prop:chains}.
Because we have assumed $\q>4$, Lemma~\ref{lemma:three_centralizers} implies that $M$ has only two distinct
two-step centralizers, and hence Lemma~\ref{lemma:2} applies.

We expand $0=[wxy[v_1xyx^{\q-2}]]=[wxy[yx^{\q-1}yx^{\q-2}]]$ by iterated application of the Jacobi identity.
Such calculation is greatly simplified if we keep in mind that all left-normed Lie brackets obtained
at the end of the procedure have $w$ as their first entry, followed by exactly three occurrences of $y$ in some places,
with $x$ filling all the remaining entries.
Because of what we have observed above, there are exactly three entries where $y$ can appear if the Lie bracket is to be nonzero,
and one of them must be one of the last two entries.
Therefore, we can disregard any complex Lie bracket in the course of the calculation as soon as it is clear that
it will not allow $y$ to appear in any of the last two entries once expanded completely.
The same applies to the other two specific entries where $y$ must appear.
All calculations in this and later proofs benefit from similar observations,
which we will not repeat in such detail.

We have $t+1=\q-2$, and according to the general theory of constituent lengths
$r+1=\q$, or $r+1=\q-2^s$ for some $0\leq s\leq e$.
In the first case we have $v=[wxyx^{\q-1}yx^{\q-2}]$, and hence
\begin{equation}\label{eq:lemma3eq1}
\begin{split}
0&=[wxy[v_{1}xyx^{\q-2}]]
\\&=
\binom{\q-2}{\q-3}[wxyx^{\q-3}[v_{1}xy]x]+\binom{\q-2}{\q-2}[wxyx^{\q-2}[v_{1}xy]]
\\&=
[wxyx^{\q-2}[v_{1}x]y]=[wxyx^{\q-1}yx^{\q-2}y]=[vxy].
\end{split}
\end{equation}
Together with $[vxx]=0$, this implies that $[vx]$ is central, a contradiction.

Now suppose that $r+1=\q-2^s$ for some $0\leq s\leq e$.
We have
\begin{equation}\label{eq:lemma3eq2}
\begin{split}
0&=[wxy[v_{1}xyx^{\q-2}]]
\\&=
\sum_{i=0}^{\q-2}\binom{\q-2}{i}[wxyx^{i}[v_{1}xy]x^{\q-2-i}]
\\&=
\sum_{i=0}^{\q-2}\binom{\q-2}{i} [wxyx^{i}[v_{1}x]yx^{\q-2-i}]
+\sum_{i=0}^{\q-2}\binom{\q-2}{i}[wxyx^{i}y[v_{1}x]x^{\q-2-i}].
\end{split}
\end{equation}
Note that the binomial coefficient $\binom{\q-2}{i}$ has the opposite parity as $i$, and hence all the terms
in the final two sums corresponding to odd values of $i$ vanish.
We distinguish two cases according to the parity of $r$.

When $r$ is odd, that is, if $1\leq s\leq e$, the only complex Lie bracket in the latter sum which may possibly not vanish
occurs when $i$ equals $r=\q-2^s-1$,
but then the corresponding binomial coefficient vanishes.
As for the former sum, note that all terms vanish except, possibly,
when the last entry $y$ is followed by exactly $2^s$ or $2^s+1$ entries $x$.
Because the corresponding binomial coefficient is nonzero only in the former case, we obtain
\begin{equation}\label{eq:lemma3eq3}
\begin{split}
0&=\sum_{i=0}^{\q-2}
\binom{\q-2}{i}[wxyx^{i}[v_{1}x]yx^{\q-2-i}]
\\&=
[wxyx^{\q-2^s-2}[v_{1}x]yx^{2^s}]
\\&=
[wxyx^{\q-2^s-1}yx^{\q-2}yx^{2^s}]
=[vxyx^{2^s}].
\end{split}
\end{equation}

Similar arguments show that in the remaining case $s=0$, which means $r=\q-2$,
all terms in the two sums vanish except, possibly, those with $i=\q-4 $ or $\q-2$ in the former sum,
and that with $i=\q-2$ in the latter.
Here $v=[wxyx^{\q-2}yx^{\q-3}]$, and hence
\begin{equation}\label{eq:lemma3eq4}
\begin{split}
0&=[wxyx^{\q-4}[v_{1}x]yxx]+[wxyx^{\q-2}[v_{1}x]y]+[wxyx^{\q-2}y[v_{1}x]]
\\&=
[vyxx]+[vxxy]+([vyxx]+[vxyx]+[vxxy])
=[vxyx].
\end{split}
\end{equation}
We have found that $[vxyx^{2^{n}}]=0$ in all cases.
But $[vxyx^i]$ spans $L_{k+i+1}$ for $0\leq i \leq q-1$ because of Lemma~\ref{lemma:2}.
Therefore, we have $L_{k+2^s+1}=\{0\}$ if $s<e$, and hence $\dim L\le k+2^s+2\le k+q+2$.

In the exceptional case $s<e$, which corresponds to $r=q-1$,
we have $v=[wxyx^{q-1}yx^{\q-3}]$.
Consider the element $[v_{1}xyx^{q-1}]$, which is nonzero, but must be centralized by either $\F y$ or $\F x$.
In the former case, the calculation
\begin{align*}
0&=
[wxyx^{q-2}[v_{1}xyx^{q-1}y]]
\\&=
[wxyx^{q-2}[v_{1}xyx^{q-1}]y]
\\&=
[wxyx^{q-1}[v_{1}xy]x^{q-2}y]
\\&=
[wxyx^{q-1}[v_{1}x]yx^{q-2}y]+[wxyx^{q-1}y[v_{1}x]x^{q-2}y]
\\&=
[vxxyx^{q-2}y]+[vyx^{q}y]
=
[vyx^{q}y],
\end{align*}
together with the fact that $[vxyx^{q}]=0$ found earlier,
implies that $[vxyx^{q-1}]$ is central, and hence $\dim L\le k+q+2$.
In the latter case, the same conclusion follows from the calculation
\[
0=[wxyx^{q-3}[v_{1}xyx^{q}]]=[wxyx^{q-3}[v_{1}xy]x^{q}]=[wxyx^{q-3}[v_{1}x]yx^{q}]=[vyx^{q}]
\]
(where $q-3\ge 0$ because $\q>4$)
together with the fact that $[vxyx^{q-2}y]=0$.
\end{proof}

Note that if the conclusion of Lemma~\ref{lemma:3} is violated, then the second constituent of $M$ is not the last,
hence $k+1\ge 2\q+q=5q\ge 20$, and so the upper bound $k+q+2$ for $\dim L$ given in Lemma~\ref{lemma:3}
does not exceed $6(k+1)/5+1$ and hence, in turn, the bound $(4k+1)/3$ of Theorem~\ref{thm:main}.

\begin{lemma}\label{lemma:4}
Suppose that the second constituent of $M$ is shorter than $\q-1$, and that the last constituent of $M$ has length $\q-1$.
Then $\dim L\le k+q+1$.
\end{lemma}

\begin{proof}
The hypotheses imply that $M$ has more than two constituents.
Also, because the length of the last constituent of $M$ is different from $\q$ and $q$,
Lemma~\ref{lemma:three_centralizers} implies that $M$ has only two distinct
two-step centralizers, and hence Lemma~\ref{lemma:2} applies.

We expand $0=[wxy[v_1xyx^{\q-2}]]=[wxy[yx^{\q-1}yx^{\q-2}]]$ as we did in the proof of  of Lemma~\ref{lemma:3}.
According to the various possible values of $r$,
the calculations involved are very similar to, but slightly simpler
than corresponding calculations in the proof of Lemma~\ref{lemma:3}.

Thus, if $r+1=\q$ we have $v=[wxyx^{\q-1}yx^{\q-2}]$, and similarly to Equation~\eqref{eq:lemma3eq1} we find
the contradiction
\begin{align*}
0&=[wxy[v_{1}xyx^{\q-2}]]
=
[wxyx^{\q-2}[v_{1}xy]]
\\&=
[wxyx^{\q-2}[v_{1}x]y]=[wxyx^{\q-1}yx^{\q-2}y]=[vy],
\end{align*}
where we have safely disregarded any Lie bracket which would not allow $y$ to appear as the last entry once expanded completely.

Now suppose that
$r+1=\q-2^s$ for some $0\leq s\leq e$.
Then Equation~\eqref{eq:lemma3eq2} holds,
and if $n>0$ the same argument employed there leads to a calculation formally identical to Equation~\eqref{eq:lemma3eq3}
except for the final result:
\[
0=[wxyx^{\q-2^s-1}yx^{\q-2}yx^{2^s}]
=[vyx^{2^s}].
\]
As in Lemma~\ref{lemma:3} it follows that $\dim L\le k+q+1$,
and the case $s=e$ needs no special treatment here.
Finally, when $s=0$ the calculation is just a little simpler than Equation~\eqref{eq:lemma3eq4}, and gives
\[
0=[wxyx^{\q-2}[v_{1}x]y]+[wxyx^{\q-2}y[v_{1}x]]
=
[vxy]+\bigl([vyx]+[vxy]\bigr)
=[vyx],
\]
a contradiction.
\end{proof}

\subsection{Second constituent of $M$ of length $\q-1$}\label{subsec:t=2q-3}

After the results of the previous subsection we may assume from now on that the second constituent of $M$ has length $\q-1$.
Consequently, according to Proposition~\ref{prop:chains},
every constituent of $M$ has length $\q$ or $\q-1$ except, possibly, for the last.
However, according to Lemma~\ref{lemma:1}, in our setting the last constituent of $M$ has length $\q-2$ or $\q-1$.
Here we prove that the former alternative leads to a contradiction.
In all the subsequent arguments, $M$ has only two distinct
two-step centralizers according to Lemma~\ref{lemma:three_centralizers}.
In particular, Lemma~\ref{lemma:2} applies.

We introduce some convenient notation.
As in the proof of Proposition~\ref{prop:chains} we set $v_1=[yx^{\q-2}]$ and $v_2=[v_1xyx^{\q-3}]$,
so that $v_{1}$ and $v_{2}$ are non zero elements of degrees $\q-1$ and $2\q-2$.
Furthermore, the relations $[v_iy]=0$ and $[v_ixx]=0$ hold, for $i=1,2$.
In fact, the rationale behind this piece of notation is that we reserve symbols such as $v$, $u$, $w$, and later $v_m$,
to denote nonzero elements in a homogeneous component immediately preceding a diamond,
possibly fake of type one, when $L$ is interpreted as a thin Lie algebra of Nottingham type
as explained in Section~\ref{sec:main_thin}.
Thus, for example, for the two elements $v$ and $u$ introduced already in Subsection~\ref{subsec:initial}
the elements $[vx]$ and $[vy]$ span a genuine diamond (the second diamond $L_k$),
while $[ux]$ spans a fake diamond of type one, and hence the relations $[uy]=0$ and $[uxx]=0$ hold.
With terminology as in~\cite{CN}, which applies to the quotient $M$ of $L$,
elements $[ux]$, $[vx]$, etc., will be at the beginning of some constituent of $M$.

\begin{lemma}\label{lemma:5}
If the second constituent of $M$ has length $\q-1$, then the last constituent has length $\q-1$ as well.
\end{lemma}

\begin{proof}
We know from Lemma~\ref{lemma:1} that the last constituent of $M$ has length $\q-2$ or $\q-1$.
Suppose for a contradiction that the former alternative holds, and recall
that $r+1$ denotes the length of the penultimate constituent of $M$.

Suppose first that $r+1=\q$.
Then we obtain a contradiction by computing
\begin{align*}
0&=[wx[v_2y]]=[wx[v_{1}xyx^{\q-3}y]]
\\&=
[wx[v_{1}xyx^{\q-3}]y]+[wxy[v_{1}xyx^{\q-3}]]
\\&=
[wx[v_{1}xy]x^{\q-3}y]+[wxyx^{\q-3}[v_{1}xy]]
\\&=
[wx[v_{1}x]yx^{\q-3}y]+[wxy[v_{1}x]x^{\q-3}y]+[wxyx^{\q-3}[v_{1}x]y]
\\&=
[wxyx^{\q-1}yx^{\q-3}y]
=[vy].
\end{align*}
Now suppose that $r+1<\q$.
Because the second constituent of $L$ has length $\q-1$,
and according to Proposition~\ref{prop:chains}, every constituent of $M$, except possibly the last,
has length at least $\q-1$.
Furthermore, if the last constituent of $M$ is shorter than $\q-1$, which is the case here, then
the penultimate constituent of $M$ has length $\q-1$.
Consequently, we have $r=\q-2$, and hence
$v=[wxyx^{\q-2}yx^{\q-3}]$.
We separately compute
\[
[wxv_{2}]=[wx[v_{1}xyx^{\q-3}]]
=
[wx[v_{1}xy]x^{\q-3}]
=
[wxy[v_{1}x]x^{\q-3}]=[vx],
\]
and
\begin{align*}
[wxyv_{2}]&=[wxy[v_{1}xyx^{\q-3}]]
\\&=
\binom{\q-3}{\q-4}[wxyx^{\q-4}[v_{1}xy]x]
+\binom{\q-3}{\q-3}[wxyx^{\q-3}[v_{1}xy]]
\\&=
[wxyx^{\q-4}[v_{1}x]yx]
+[wxyx^{\q-3}[v_{1}x]y]
\\&=
[vyx]+[vxy].
\end{align*}
It follows that
$0=[wx[v_{2}y]]=[wxv_{2}y]+[wxyv_{2}]=[vyx]$, again a contradiction,
under our blanket assumption that $\dim L\ge k+3$.
\end{proof}

\subsection{The structure of $M$ past the second diamond}\label{subsec:t=2q-2}

We may now assume that both the second and the last constituent of $M$ have length $\q-1$.
We know from Proposition~\ref{prop:chains} that every constituent of $M$ has length $\q$ or $\q-1$.
Our remaining goals are to prove that all constituents of $M$ except the second and the last have length $\q$,
and that the total number of constituents of $M$ equals a power of two.
Reversing the logical order, we prove the latter assertion first, assuming the former.
This is more convenient because it correctly shapes our notation
to describe the initial structure of genuine examples of large dimension, while our proof of the former assertion
in Subsection~\ref{subsec:final} will only produce a contradiction in the case to be excluded.
A further reason for postponing Subsection~\ref{subsec:final} is that it is the most complex part of the proof.

Because we will work with a considerable portion of $L$ past the second diamond, we conveniently switch
from the terminology of graded Lie algebras of maximal class to that of thin Lie algebras of Nottingham type.
Thus, assume that $k=n \q-1$ for some $n\geq 3$, and that $L$, after the fake second diamond in degree $q$, has a fake
diamond of type one in each degree $m \q-1$, for $1<m<n$.
Recall from Section~\ref{sec:main_thin} that a diamond of type one is a one-dimensional component $\F [wx]$ of $L$ such that
$[wy]=0$ and $[wxx]=0$;
in other words, it is a one-dimensional component centralized by $x$ lying between two components centralized by $y$.
Our goal is to show that $n$ is a power of two unless $L$ has bounded dimension.
To achieve this we need to gather a certain amount of information on $L$ past the diamond $L_k$.
In particular, we show that the genuine diamond $L_k$ is followed by a certain number of fake diamonds of type one
occurring at regular intervals.

We have already introduced elements $v_1=[yx^{\q-2}]$ and $v_2=[v_1xyx^{\q-3}]$, of degrees $\q-1$ and $2\q-2$.
We extend these definitions by recursively setting
\[
v_{m}=[v_{m-1}xyx^{\q-2}] \quad\text{for}\quad 2<m\leq n,
\]
so that $v_{m}$ is a nonzero element of  degree $m\q-2$ for $m>1$, and if $m<n$ it satisfies
\[
[v_{m}y]=0,\qquad [v_{m}xx]=0, \quad\text{and}\quad [v_{m}xyx^{j}y]=0 \quad\text{for}\quad 0\le j<\q-2.
\]
It is our goal, in the most part of this subsection, to prove that this periodic structure
of $L$ is only slightly perturbed when passing through the diamond $L_k$, and continues for a while after it.
Note that the first equation, $[v_{m}y]=0$ is the same as
$[v_{m-1}xyx^{\q-2}y]=0$ for $m>1$, and hence could be included in the third set of relations by extending the range of $j$,
but such a formulation would be less convenient for our purposes.
Note also that $v_n$ equals, up to a scalar, the element denoted by $v$ in Subsection~\ref{subsec:initial}.
Hence $[v_ny]$ is nonzero here and, in fact together with $[v_nx]$ it spans the second genuine diamond $L_k$.
According to Lemma~\ref{lemma:2} we have
\[
[v_{n}xx]=0, \quad\text{and}\quad [v_{n}xyx^{j}y]=0 \quad\text{for}\quad 0\le j<q-1,
\]
which have been crucial facts in most arguments so far.
Under the present more special assumptions, a different calculation allows us to extend the latter assertion to $0\le j<\q-2$,
in Lemma~\ref{lemma:6}.
Generally speaking, here and in similar situations later, a proof for roughly the first half of the range for $j$
may be achieved by using only the first constituent of $M$, as in the proof of Lemma~\ref{lemma:2},
but to deal with the whole range one needs to exploit the first two constituents of $M$.
More precisely, one uses the fact that the second constituent attains the maximum length $\q-1$
allowed by the general theory, as in the proof of Proposition~\ref{prop:chains}.

Before proceeding further we note that, because $[v_nxx]=0$ (and $[v_nyy]=0$), the quotient $L/L^{k+2}$ is $\Z^2$-graded
by assigning independent degrees to $x$ and $y$.
This fact will extend to each graded quotient of $L$ which we will implicitly consider in
all the calculations of the rest of the proof.
As observed at the beginning of the proof of Proposition~\ref{prop:chains},
this implies that, at any stage where the structure of $L/L^{s+1}$ has been determined and $\dim L_s=1$,
in any nonvanishing long Lie bracket of length $s$
the symbols $x$ and $y$ must appear the expected number of times, and at specific places
(allowing for both $[v_nxy]$ and $[v_nyx]$ in degree $k+1$).

\begin{lemma}\label{lemma:6}
We have
$[v_{n}xyx^{j}y]=0$
for $0\le j<\q-2$.
\end{lemma}

\begin{proof}
We prove the result by induction on $j$ using the equation $[v_{1}xyx^{j}y]=0$,
which holds for $0\le j<\q-2$.
The induction step follows from the calculation
\begin{align*}
0&=[v_{n-1}x[v_{1}xyx^{j}y]]
\\&=
[v_{n-1}x[v_{1}xyx^{j}]y]+
[v_{n-1}xy[v_{1}xyx^{j}]]
\\&=
[v_{n-1}x[v_{1}xy]x^{j}y]+[v_{n-1}xyx^{j}[v_{1}xy]]
\\&=
[v_{n-1}x[v_{1}x]yx^{j}y]+[v_{n-1}xy[v_{1}x]x^{j}y]
+[v_{n-1}xyx^{j}[v_{1}x]y]
\\&=
[v_{n}xyx^{j}y]
+\bigl([v_{n}yxx^{j}y]+[v_{n}xyx^{j}y]\bigr)
+\bigl([v_{n}yxx^{j}y]+[v_{n}xyx^{j}y]\bigr)
\\&=
[v_{n}xyx^{j}y].
\end{align*}
As in previous occasions, we have repeatedly used the fact that any long Lie bracket involved
ending in $x$ vanishes, because then $y$ occurs too many times before the end,
against the induction hypothesis.
\end{proof}

As a first step towards our conclusion that $n$ is a power of two for $\dim L$ large, we prove that $n$ is even.
Generally speaking, in this subsection we need to use the initial structure of $M$ up to its second constituent.
Therefore, the adjoint action on $L$ of elements near the end of the first or second constituent of $M$
will be important.
To begin with, short calculations show that
\[
[v_{m}x[v_{1}x]]=[v_{m+1}x] \quad \textrm{and} \quad [v_{m}x[v_{2}x]]=[v_{m+1}yx^{\q-1}]
\]
for $2\le m<n$.
In particular, the former equation yields the convenient formula $[v_mx]=[v_2x[v_1x]^{m-2}]$ by induction,
and according to the latter equation we have $[v_{m}x[v_{2}x]]=0$ for $m<n-1$.
If $n$ is odd, then the calculation
\begin{align*}
0&=
[v_{(n+1)/2}x[v_{(n+1)/2}x]]
\\&=
[v_{(n+1)/2}x[v_{2}x[v_{1}x]^{(n-3)/2}]]
\\&=
[v_{(n+1)/2}x[v_{1}x]^{(n-3)/2}[v_{2}x]]
\\&=
[v_{n-1}x[v_{2}x]]=[v_{n}yx^{\q-1}]
\end{align*}
shows that $[v_{n}xyx^{\q-2}]=0$.
Together with Lemma~\ref{lemma:6} this implies that $\dim L\le k+\q$.

Assuming $\dim L$ larger than that, we have proved that $n$ must be even.
To proceed further we need to prove that all homogeneous components of $L$ up to a certain degree
are centralized by either $x$ or $y$ (whence they are one-dimensional) and that the genuine second diamond $L_k$
is followed by a certain number of fake diamonds of type one in degrees of the form $m\q-2$.
The proof will be by induction on $m$, but the initial step requires special treatment.
Set
\[
v_{n+1}=[v_{n}xyx^{\q-3}],
\]
so that $v_{n+1}$ is a nonzero element of  degree $(n+1)\q-3$.
Note that the difference in degrees between $v_n$ and $v_{n+1}$ is only $\q-1$,
hence one less than between earlier consecutive $v_m$ (with the exception of that between $v_1$ and $v_2$).
Recall from the end of Subsection~\ref{subsec:initial} that $[v_{n}yx]=(\mu^{-1}-1)[v_{n}xy]$, where $\mu$
is the type of the second genuine diamond $L_k$, and so $\mu^{-1}-1=\mu^{-1}+1\in\F\setminus\{0\}$.
Note that we do not need to require that $\dim L$ is large enough in the following two lemmas,
as the claimed equations hold trivially otherwise.

\begin{lemma}\label{lemma:7}
We have
\[
[v_{n+1}y]=0,\qquad [v_{n+1}xx]=0, \quad\text{and}\quad [v_{n+1}xyx^{j}y]=0 \quad\text{for}\quad 0\le j<\q-2.
\]
\end{lemma}

\begin{proof}
The equation $[v_{n+1}y]=0$ was included in Lemma~\ref{lemma:6} as the case $j=\q-3$.
The equation $[v_{n+1}xx]=0$ follows from
\[
0=[v_{n-1}x[v_{2}xx]]=[v_{n-1}x[v_{2}x]x]=[v_{n}yx^{\q-1}x]=(\mu^{-1}-1)[v_{n+1}xx].
\]

It remains to prove that $[v_{n+1}xyx^{j}y]=0$ for $0\le j<\q-2$, by induction on $j$.
A calculation analogous to that in the proof of Lemma~\ref{lemma:6},
with $n-1$ replaced by $n$, would be inconclusive here.
However, the slight variation
\begin{align*}
0&=[v_{n}[v_{1}xyx^{j}y]]=[v_{n}[v_{1}xyx^{j}]y]+[v_{n}y[v_{1}xyx^{j}]]
\\&=
[v_{n}[v_{1}xy]x^{j}y]+j[v_{n}x[v_{1}xy]x^{j-1}y]+[v_{n}yx^{j}[v_{1}xy]]
\\&=
\bigl([v_{n}[v_{1}x]yx^{j}y]+[v_{n}y[v_{1}x]x^{j}y]\bigr)
\\&\quad
+j\bigl([v_{n}x[v_{1}x]yx^{j-1}y]+[v_{n}xy[v_{1}x]x^{j-1}y]\bigr)
+[v_{n}yx^{j}[v_{1}x]y]
\\&=
\bigl(\mu^{-1}[v_{n+1}xyx^{j}y]+(\mu^{-1}-1)[v_{n+1}xyx^{j}y]\bigr)
\\&\quad
+j[v_{n+1}xyx^{j}y]
+(\mu^{-1}-1)[v_{n+1}xyx^{j}y]
\\&=
(\mu^{-1}+j)[v_{n+1}xyx^{j}y]
\end{align*}
proves the induction step for $j$ odd (because $\mu^{-1}\neq 1$).
But for even values of $j$ the much shorter calculation
\[
0=[v_{n+1}[yx^{j+1}y]]=(j+1)[v_{n+1}xyx^jy]
\]
gives the desired conclusion.
\end{proof}

We extend the definition of the elements $v_m$ past $v_{n+1}$ by recursively setting
\[
v_{m}=[v_{m-1}xyx^{\q-2}], \quad\text{for}\quad n+1<m\le 3n/2,
\]
so that $v_{m}$ acquires formal degree $m\q-3$ for $m>n$,
We now prove that for $m<3n/2$ these elements $v_m$ satisfy analogous equations as those with $m<n$.

\begin{lemma}\label{lemma:8}
For $n<m<3n/2$ we have
\[
[v_{m}y]=0,\qquad [v_{m}xx]=0, \quad\text{and}\quad [v_{m}xyx^{j}y]=0 \quad\text{for}\quad 0\le j<\q-2.
\]
\end{lemma}

\begin{proof}
The proof proceeds by induction on $m$, the case $m=n+1$ being Lemma~\ref{lemma:7}.
Thus, assume that the conclusion holds for $m<n+\ell$, where $2\leq\ell< n/2$.
We now prove that it holds for $m=n+\ell$ as well.

The first equation, $[v_{n+\ell}y]=0$, is actually the hardest to prove,
as its proof possibly relies on the global structure of $L$, that is, in all previous degrees.
In fact, that is the equation which determines the type of the diamond spanned by $[v_mx]$ and $[v_mx]$,
and is the obstacle to extending the result to $m=3n/2$, see our comments after this proof.
Hence we postpone a proof of $[v_{n+\ell}y]=0$ and show first that the other equations follow from it
together with local calculations.
The equation
$[v_{n+\ell}xx]=0$ follows at once from $0=[v_{n+\ell-1}x[v_{1}xx]]=[v_{n+\ell-1}x[v_{1}x]x]$.
The last set of equations
can be proved by induction on $j$ using the relation $[v_{1}xyx^{j}y]=0$, for $0\le j<\q-2$,
in a similar way as in the proof of Lemma~\ref{lemma:6},
the induction step being
\begin{align*}
0&=[v_{n+\ell-1}x[v_{1}xyx^{j}y]]
\\&=
[v_{n+\ell-1}x[v_{1}xyx^{j}]y]+
[v_{n+\ell-1}xy[v_{1}xyx^{j}]]
\\&=
[v_{n+\ell-1}x[v_{1}xy]x^{j}y]+[v_{n+\ell-1}xyx^{j}[v_{1}xy]]
\\&=
[v_{n+\ell-1}x[v_{1}x]yx^{j}y]+[v_{n+\ell-1}xy[v_{1}x]x^{j}y]
+[v_{n+\ell-1}xyx^{j}[v_{1}x]y]
\\&=
3[v_{n+\ell}xyx^{j}y].
\end{align*}

It remains to prove that $[v_{n+\ell}y]=0$,
but to facilitate our calculations, we first collect some information on the adjoint action
of $[v_1x]=[yx^{\q-1}]$ and $v_2=[v_1xyx^{\q-3}]$.
For the former one easily finds that
\begin{equation}\label{eq:v_m}
[v_{m}xyx^{j}[v_1x]]=[v_{m+1}xyx^{j}] \quad\text{for}\quad 0\le j\le\q-1
\end{equation}
for $2\le m<n-1$ or $n<m<n+\ell-1$.
Note that this equation reads $[v_{m+1}[v_1x]]=v_{m+2}$ and $[v_{m+1}x[v_1x]]=[v_{m+2}x]$ when $j=\q-2$ or $\q-1$,
these two formulas holding when $m+1=2$ as well.
However, we have
\begin{equation}\label{eq:v_n-1}
[v_{n-1}xyx^{j}[v_1x]]=\mu^{-1}[v_{n}xyx^{j}] \quad\text{for}\quad 0\le j<\q-1,
\end{equation}
(which reads $[v_n[v_1x]]=\mu^{-1}[v_{n+1}x]$ when $j=\q-2$),
$[v_nx[v_1x]]=0$,
$[v_ny[v_1x]]=(\mu^{-1}+1)[v_{n+1}xy]$,
and
\begin{equation*}
[v_{n}xyx^{j}[v_1x]]=[v_{n+1}xyx^{j+1}] \quad\text{for}\quad 0\le j<\q-1.
\end{equation*}
The action of
$v_2=[v_1xyx^{\q-3}]$
is slightly more complex to describe in general, and we limit ourselves to the information which we need here.
For $2\le m<n+\ell-2$ with $m\neq n-1, n$ we have
\begin{equation}\label{eq:v_mxv_2}
[v_{m}xv_{2}]
=[v_{m}x[v_{1}xy]x^{\q-3}]
=[v_{m}x[v_{1}x]yx^{\q-3}]+[v_{m}xy[v_{1}x]x^{\q-3}]
=0,
\end{equation}
because the two summands cancel out,
while
\begin{equation}\label{eq:v_n-1xv_2}
\begin{split}
[v_{n-1}xv_{2}]
&=
[v_{n-1}x[v_{1}x]yx^{\q-3}]+[v_{n-1}xy[v_{1}x]x^{\q-3}]
\\&=
[v_{0}xyx^{\q-3}]
+\bigl([v_{0}yxx^{\q-3}]
+[v_{0}xyx^{\q-3}]\bigr)
=(1+\mu^{-1})v_{n+1}
\end{split}
\end{equation}
and
\begin{equation}\label{eq:v_nxv_2}
[v_{n}xv_{2}]=[v_{n}x[v_{1}xy]x^{\q-3}]=[v_{n+1}xyx^{\q-2}]=v_{n+2}.
\end{equation}
We also have
$[v_mxyv_2]=0$ and
$[v_mxyxv_2]=0$
for $2\le m<n+\ell-2$,
with $m\neq n-2,n-1$ in the latter.

With these formulas at hand we are ready to complete the induction step by proving that $[v_{n+\ell}y]=0$.
For $0<s<n-\ell$ (whence $n-s\ge 2$) we compute
\begin{align*}
0&=
[v_{n-s}x[v_{\ell+s}y]]
\\&=
[v_{n-s}xv_{\ell+s}y]+[v_{n-s}xyv_{\ell+s}]
\\&=
[v_{n-s}x[v_{2}[v_{1}x]^{\ell+s-2}]y]
+[v_{n-s}xy[v_{2}[v_{1}x]^{\ell+s-2}]]
\\&=
\sum_{i=0}^{\ell+s-2}\binom{\ell+s-2}{i}
[v_{n-s}x[v_{1}x]^{i}v_{2}[v_{1}x]^{\ell+s-2-i}y]
\\&\quad+
\sum_{i=0}^{\ell+s-2}\binom{\ell+s-2}{i}
[v_{n-s}xy[v_{1}x]^{i}v_{2}[v_{1}x]^{\ell+s-2-i}]
\\&=
(\mu^{-1}+1)
\binom{\ell+s-2}{s-1}
[v_{n+1}[v_{1}x]^{\ell-1}y]
+\binom{\ell+s-2}{s}
[v_{n+2}[v_{1}x]^{\ell-2}y]
\\&\quad+
\mu^{-1}
[v_{n-s}xy[v_{1}x]^{\ell+s-2}v_2]
\\&=
\left((\mu^{-1}+1)\binom{\ell+s-2}{s-1}
+\binom{\ell+s-2}{s}
+\mu^{-1}\right)[v_{n+\ell}y],
\end{align*}
where most of the summands vanish because of terms centralized by $v_2$.
Now it suffices to find a value of $s$, depending on $\ell$, such that the final coefficient is nonzero.
Taking $s=\ell$, which is allowed because $\ell<n/2$, the coefficient of $[v_{n+\ell}y]$ equals
$\binom{2\ell-2}{\ell}+\mu^{-1}$,
because
$\binom{2\ell-2}{\ell-1}=\binom{2\ell-3}{\ell-2}+\binom{2\ell-3}{\ell-1}=2\binom{2\ell-3}{\ell-1}$
is even for any $\ell>1$.
Recalling that $\mu^{-1}+1$ is a nonzero element of $\F$,
that coefficient does not belong to the prime field $\F_2$, and hence is not zero, if $\mu^{-1}\neq 0$.
In the remaining case where $\mu^{-1}=0$ (which means that the diamond $L_k$ has type $\infty$)
the coefficient of $[v_{n+\ell}y]$ equals
$\binom{\ell+s-2}{s-1}+\binom{\ell+s-2}{s}=\binom{\ell+s-1}{s}$
in $\F$, which is nonzero, for example, when $s$ equals the $2$-part of $\ell$.
\end{proof}

The equations proved in Lemma~\ref{lemma:8} say that, after the genuine diamond in degree $k=n\q-1$, if $\dim L$ is large enough,
$L$ has diamonds of type one in each degree $m\q-2$ with $n<m<3n/2$.
This cannot be further extended in general, as Young's construction from~\cite{Young:thesis}
allows one to construct (infinite-dimensional) examples where
$L$ has a genuine diamond in degree $(3n/2)\q-2$.
Hence $[v_{3n/2}y]$ does not vanish in those algebras.

We are finally able to bound $\dim L$ if $n$ is not a power of two.
Assuming $n$ even as we may, we set $t=(n+n_2)/2$, where $n_2$ denotes the $2$-part of $n$.
Suppose that $n$ is not a power of two.
Then $n\ge 3n_2$, hence $2t\le 4n/3<3n/2$,
and so Lemma~\ref{lemma:8} certainly holds for $n<m\le 2t$.
In particular, we have $[v_{2t}y]=0$.
Using Equations~\eqref{eq:v_mxv_2}, \eqref{eq:v_n-1xv_2} and~\eqref{eq:v_nxv_2} we compute
\begin{align*}
0&=
[v_{t}x[v_{t}x]]
=
[v_{t}xv_{t}x]
\\&=
[v_{t}x[v_{2}[v_{1}x]^{t-2}]x]
\\&=
\sum_{i=0}^{t-2}\binom{t-2}{i}
[v_{t}x[v_{1}x]^{i}v_{2}[v_{1}x]^{t-2-i}x]
\\&=
\binom{t-2}{n-1-t}
[v_{n-1}v_{2}[v_{1}x]^{2t-n-1}x]
+
\binom{t-2}{n-t}
[v_nv_{2}[v_{1}x]^{2t-n-2}x]
\\&=
(\mu^{-1}+1)\binom{t-2}{2t-n-1}[v_{2t}x]+\binom{t-2}{2t-n-2}[v_{2t}x].
\end{align*}
Now Lucas' theorem implies that
$
\binom{t-2}{2t-n-1}
=
\binom{(n-n_2)/2+(n_2-2)}{n_2-1}
$
is even and
$
\binom{t-2}{2t-n-2}
=
\binom{(n-n_2)/2+(n_2-2)}{n_2-2}
$
is odd, noting that $(n-n_2)/2$ is a multiple of $n_2$.
We conclude that $[v_{2t}x]=0$,
and hence $L$ has bounded dimension.
More precisely, in the worst case where $v_{2t}\neq 0$ we have $\dim L=(n+n_2)\q-1=k+(k+1)_2$.
Because $k+1\ge 3(k+1)_2$ if $k+1$ is not a power of two, it follows that $\dim L\le (4k+1)/3$.
Consequently, if $\dim L$ exceeds $(4k+1)/3$ as required in Theorem~\ref{thm:main} then $n$,
and with it $k+1$, must be a power of two.
Note that this is the calculation in highest degree that we have considered so far.

Modulo the part of the proof which we have postponed to Subsection~\ref{subsec:final},
this completes the proof of Theorem~\ref{thm:main}, and the equivalent Theorem~\ref{thm:main_thin},
and also the supplementary assertions~(4) and~(5) of Theorem~\ref{thm:main_thin_more}.
Now we can use the information that $n$ is a power of two larger than two
to prove assertion~(6) of Theorem~\ref{thm:main_thin_more}.
In fact, we only need to know that $n$ is a multiple of four,
which is granted as soon as $\dim L>k+2\q$ according to the previous paragraph.
Thus, if $\dim L>k+2\q$ then $n/2$ is even, and the calculation
\begin{align*}
0&=
[v_{n/2+1}x[v_{n/2+1}x]]
\\&=
[v_{n/2+1}x[v_2x[v_1x]^{n/2-1}]]
\\&=
(n/2-1)[v_{n-1}x[v_2x][v_1x]]
+[v_nx[v_2x]]
\\&=
[v_{n-1}xv_2x[v_1x]]
+[v_nxv_2x]
\\&=
(\mu^{-1}+1)[v_{n+2}x]
+[v_{n+2}x]
=
\mu^{-1}[v_{n+2}x]
\end{align*}
implies that $\mu=\infty$, as desired.

\subsection{Eliminating irregular configurations}\label{subsec:final}

It remains to prove that
if $\dim L$ is large enough then
all constituents of $M$ apart from the second and the last have has length $\q$.
Suppose that, on the contrary, $M$ has at least three constituents of length $\q-1$.
Thus, the sequence of constituent lengths of $M$ reads
\[\boxed{
\q,\q-1,\ldots, \q-1, \q^{\ r-2},\q-1,
}\]
where $\q^{\ r-2}$ indicates a sequence of $r-2$ consecutive occurrences of $\q$, for some $r\geq 2$.
We will gather more and more of the missing information in the sequence of constituent lengths of $M$,
both from the beginning and the end, until we manage to obtain a contradiction in all cases, for $\dim L$ large enough.

Fix a nonzero element $v_0$ of $L_{k-1}$.
This plays the role which was of $v_n$ in the setting of Subsection~\ref{subsec:t=2q-2}.
In particular, Lemma~\ref{lemma:6} remains valid here, but reads
$[v_{0}xyx^{j}y]=0$
for $0\le j<\q-2$.
Consequently, $[v_{0}xyx^{j}]$ spans $L_{k+j}$ for $0\le j\le\q-2$.
For $0<i<r$, define inductively $v_{-i}$
as the unique homogeneous element of $L$ such that $[v_{-i}xyx^{\q-2}]=v_{-i+1}$.
Finally, let $v_{-r}$ be the homogeneous element such that $[v_{-r}xyx^{\q-3}]=v_{-r+1}$.
Hence $[v_{-r}x]$ is at the beginning of the penultimate constituent of $M$ of length $\q-1$.
Retaining the notation $v_1=[yx^{\q-2}]$ and $v_2=[v_1xyx^{\q-3}]$ from the previous subsections,
Equation~\eqref{eq:v_m} holds here for $-r<m<-1$,
Equation~\eqref{eq:v_n-1} holds if we read $-1$ and $0$ in place of the subscripts $n-1$ and $n$,
and the adjoint action of $[v_1x]$ on elements close to $v_{-r}$ admits a similar description.
In particular, we have $[v_{-r}[v_1x]]=[v_{-r+1}x]$ and $[v_{-r}x[v_1x]]=0$, but
$[v_{m}[v_1x]]=v_{m+1}$ for $-r<m<0$.

Because
$[v_{-r}x[v_2x]]=[v_{-r}x[v_1xy]x^{\q-2}]
=[v_{-r}xyx^{\q-3}yx^{\q-1}]=[v_{-r+2}x]$,
we have
\[
[v_{-r}x[v_2x[v_1x]^{r-2}]]
=[v_{-r}x[v_2x][v_1x]^{r-2}]
=[v_0x]\neq 0,
\]
and hence $[v_2x[v_1x]^{r-2}]\neq 0$.
This means that the second constituent of $M$ is followed by at least $r-2$ constituents of length $\q$.
If it is followed by that number and no more, that is, if the sequence of constituent lengths of $M$ begins with
$\q,\q-1,\q^{r-2},\q-1,\ldots$, then we have $[v_2x[v_1x]^{r-1}]=0$ or, equivalently, $[v_2[v_1x]^{r-1}x]=0$.
Equation~\eqref{eq:v_mxv_2} now takes the form
\begin{equation}\label{eq:v_mxv_2-bis}
[v_mxv_2]=0\quad\text{for}\quad -r<m<-1,
\end{equation}
while Equation~\eqref{eq:v_n-1xv_2} reads
$[v_{-1}xv_{2}]
=(\mu^{-1}+1)[v_0xyx^{\q-3}]$.
We also have
\begin{equation}\label{eq:v_-rxv_2}
[v_{-r}xv_2]
=[v_{-r}x[v_1xy]x^{\q-3}]
=[v_{-r}xy[v_1x]x^{\q-3}]
=[v_{-r+1}xyx^{\q-2}]
=v_{-r+2}
\end{equation}
and
\begin{equation}\label{eq:v_-rv_2}
\begin{split}
[v_{-r}v_2]
&=
[v_{-r}[v_1xy]x^{\q-3}]
+(\q-3)[v_{-r}x[v_1xy]x^{\q-4}]
\\&=
[v_{-r}[v_1x]yx^{\q-3}]
+[v_{-r}xy[v_1x]x^{\q-4}]
\\&=
[v_{-r+1}xyx^{\q-3}]
+[v_{-r+1}xyx^{\q-3}]
=0.
\end{split}
\end{equation}
After computing
\begin{align*}
0&=[v_{-r}[v_2[v_1x]^{r-1}x]]
\\&=
[v_{-r}[v_2[v_1x]^{r-1}]x]
+[v_{-r}x[v_2[v_1x]^{r-1}]]
\\&=
[v_{-r}[v_1x]^{r-1}v_2x]
+[v_{-r}xv_2[v_1x]^{r-1}]
\\&=
[v_{-1}xv_2x]
+[v_{-r+2}[v_1x]^{r-1}]
\\&=
(1+\mu^{-1})[v_0xyx^{\q-2}]
+\mu^{-1}[v_{0}xyx^{\q-2}]
\\&
=[v_0xyx^{\q-2}],
\end{align*}
we conclude that $\dim L\le k+\q$ in this case.

Hence we may assume that the second constituent of $M$ is followed by at least $r-1$ constituents of length $\q$,
and so the sequence of constituent lengths of $M$ reads
\[\boxed{
\q,\q-1, \q^{r-1},\ldots,\q-1,\q^{r-2},\q-1
}\]
Now we prove that $\dim L$ is bounded unless $r$ is even.
Here we have $[v_2[v_1x]^{r-1}x]=[v_2x[v_1x]^{r-1}]\neq 0$, but we can use $[v_2[v_1x]^{r-1}y]=0$ instead.
Note that
$[v_{-r}xyv_2]=0$ and $[v_mxyxv_2]=0$ for $-r<m<-2$
(as in the formulas right after Equation~\eqref{eq:v_nxv_2}), while
assuming $r>2$ as we may we have
\begin{equation}\label{eq:v_-1xyxv_2}
\begin{split}
[v_{-2}xyxv_2]
&=
[v_{-2}xyx[v_1xyx^{\q-3}]]
\\&=
[v_{-2}xyx^{\q-2}[v_1xy]]
=
[v_{-1}[v_1xy]]
\\&=
[v_{-1}[v_1x]y]
=[v_0y],
\end{split}
\end{equation}
and
\begin{align*}
[v_{-1}xyxv_2]
&=
[v_{-1}xyx[v_1xyx^{\q-3}]]
\\&=
[v_{-1}xyx^{\q-2}[v_1xy]]
=
[v_0[v_1xy]]
\\&=
[v_0[v_1x]y]]
+[v_0y[v_1x]]
\\&=
\mu^{-1}[v_0xyx^{\q-2}y]
+(\mu^{-1}+1)[v_0xyx^{\q-2}y]
=[v_0xyx^{\q-2}y].
\end{align*}
In these last two calculations, like elsewhere, we have used the fact that any of the resulting Lie brackets ending in $x$ vanishes
because it is preceded by too many occurrences of $y$.
Using these equations we compute
\begin{align*}
0&=[v_{-r}x[v_2[v_1x]^{r-1}y]]
\\&=
[v_{-r}x[v_2[v_1x]^{r-1}]y]
+[v_{-r}xy[v_2[v_1x]^{r-1}]]
\\&=
[v_{-r}xv_2[v_1x]^{r-1}y]
+(r-1)[v_{-r}xy[v_1x]^{r-2}v_2[v_1x]]
+[v_{-r}xy[v_1x]^{r-1}v_2]
\\&=
[v_{-r+2}[v_1x]^{r-1}y]
+(r-1)[v_{-2}xyxv_2[v_1x]]
+[v_{-1}xyxv_2]
\\&=
[v_0[v_1x]y]
+(r-1)[v_0y[v_1x]]
+[v_0xyx^{\q-2}y]
\\&=
\mu^{-1}[v_0xyx^{\q-2}y]
+(r-1)(\mu^{-1}+1)[v_0xyx^{\q-2}y]
+[v_0xyx^{\q-2}y]
\\&=
r(\mu^{-1}+1)[v_0xyx^{\q-2}y].
\end{align*}
Because $[v_0xyx^{\q-1}]=0$ as well
(being $[v_{n+1}xx]=0$ in the notation of Lemma~\ref{lemma:7}),
if $r$ is odd we conclude that $\dim L\le k+\q+1$.

Thus, we may now assume that $r$ is even.
This piece of information allows us to show that the penultimate constituent of $M$ of length $\q-1$
cannot be immediately preceded by
another constituent of length $\q-1$.
That is, we show that if the sequence of constituent lengths of $M$ reads
\[\boxed{
\q,\q-1, \q^{r-1},\ldots,\q-1,\q-1,\q^{r-2},\q-1,
}\]
then $\dim L$ is bounded.
In fact, in that case we extend our definition of the elements $v_m$ by letting $v_{-r-1}$ be
the unique homogeneous element such that
$[v_{-r-1}xyx^{\q-3}]=v_{-r}$,
so that $[v_{-r-1}x]$ is at the beginning of the third constituent of $M$ of length $\q-1$ counting from the end.
Then
\begin{align*}
[v_{-r-1}v_2]
&=
[v_{-r-1}[v_1xyx^{\q-3}]]
\\&=
[v_{-r-1}[v_1x]yx^{\q-3}]
+(\q-3)[v_{-r-1}xy[v_1x]x^{\q-4}]
\\&=
v_{-r+1}+v_{-r+1}=0
\end{align*}
and, because $[v_2[v_1x]^{r-1}y]=0$, the calculation
\begin{align*}
0&=[v_{-r-1}[v_2[v_1x]^{r-1}y]]
\\&=
[v_{-r-1}[v_2[v_1x]^{r-1}]y]
\\&=
(r-1)[v_{-r-1}[v_1x]v_2[v_1x]^{r-2}y]
\\&=
(r-1)[v_0y]
\end{align*}
yields a contradiction.

To recap, we may assume that the sequence of constituent lengths of $M$ reads
\[\boxed{
\q,\q-1, \q^{r-1},\ldots,\q-1,\q^{r-2},\q-1,
}\]
with $r$ even, and that the penultimate constituent of length $\q-1$ is not immediately preceded
by another constituent of the same length.
Now the initial section $\q,\q-1,\q^{r-1}$ of this sequence can be immediately followed by a constituent of length $\q-1$,
or by a further constituent of length $\q$.
We will be able to bound $\dim L$ in both cases, but to do that it will be necessary to go a little deeper after
the second diamond $L_k$.
To this purpose, we note that the proofs of Lemmas~\ref{lemma:7} and~\ref{lemma:8} depend only
on an initial and a final portion of the structure of $L$ up to the second diamond
(of size depending on $m$ in the latter result), and so remain valid to some extent in the present more general setting.
More precisely, resuming our notation $v_n$ for $v_0$ from Subsection~\ref{subsec:t=2q-2},
whence $v_{n+1}$ denotes $[v_{0}xyx^{\q-3}]$ (not to be confused with $v_{1}=[yx^{\q-2}]$) and then
$v_{m}=[v_{m-1}xyx^{\q-2}]$ recursively,
then Lemma~\ref{lemma:7} remains valid, and so do the assertions of Lemma~\ref{lemma:8} for $n<m<n+r/2+1$.
In particular, if $r>2$ we have $[v_{n+2}y]=0$ and $[v_{n+2}xx]=0$.
The case $r=2$ will require special treatment in a couple of places.

Suppose first that the sequence of constituent lengths of $M$ reads
\[\boxed{
\q,\q-1, \q^{r-1},\q-1,\ldots,\q-1,\q^{r-2},\q-1,
}\]
with $r$ even, which we allow to include the case $\q,\q-1, \q^{r-1},\q-1,\q^{r-2},\q-1$.
We will use the fact that $[v_2[v_1x]^rx]=0$, and expand
\[
0=[v_{-r}[v_2[v_1x]^rx]]
=
[v_{-r}[v_2[v_1x]^r]x]
+[v_{-r}x[v_2[v_1x]^r]].
\]
We compute the two summands separately.
Using Equations~\eqref{eq:v_-rv_2}, \eqref{eq:v_mxv_2-bis} and~\eqref{eq:v_nxv_2}, and the fact that $r$ is even, we obtain
\[
[v_{-r}[v_2[v_1x]^r]x]
=r[v_{-r}[v_1x]^{r-1}v_2[v_1x]x]
+[v_{-r}[v_1x]^rv_2x]
=
[v_0xv_2x]
=
[v_{n+2}x].
\]
Because of $[v_{-r}[v_1x]]=0$, Equation~\eqref{eq:v_-rxv_2} and $[v_0[v_1x]]=\mu^{-1}[v_{n+1}x]$
(a special case of Equation~\eqref{eq:v_n-1}), we have
\[
[v_{-r}x[v_2[v_1x]^r]]
=[v_{-r}xv_2[v_1x]^r]
=[v_{-r+2}[v_1x]^r]
=[v_{0}[v_1x]^2]
=\mu^{-1}[v_{n+2}x].
\]
Putting the pieces together we find that
$0=(1+\mu^{-1})[v_{n+2}x]$.
Because $[v_{n+2}y]=0$ if $r>2$ we conclude that $\dim L\le k+2\q$ in that case.

When $r=2$ we need a different argument.
Recall that the penultimate constituent of length $\q-1$ is not immediately preceded
by another constituent of the same length.
Hence when $r=2$ the sequence of constituent lengths of $M$ reads
$\q,\q-1, \q,\q-1,\ldots,\q,\q-1,\q-1$.
It is therefore natural to let $v_{-3}$ be the unique homogeneous element such that
$[v_{-3}xyx^{\q-2}]=v_{-2}$.
We would need more information to define an element $v_{-4}$ appropriately,
but here we only need an element which lies roughly half way between $v_{-4}$ and $v_{-3}$.
Thus, let $t$ be the unique homogeneous element such that
$[tx^{q-1}]=v_{-3}$, and expand
\[
0=
[t[v_2[v_1x]^2x]]
=[t[v_2x[v_1x]^2]]
=[t[v_2x][v_1x]^2]+[t[v_1x]^2[v_2x]].
\]
We find that
\[
[t[v_2x]]
=[t[v_1xyx^{\q-2}]]
=\binom{\q-2}{q-1}[tx^{q-1}[v_1xy]x^{q-1}]
=0,
\]
noting that $[v_1xy]$ centralizes almost all homogeneous elements involved,
and
\begin{align*}
[t[v_1x]^2[v_2x]]
&=
[v_{-2}xyx^{q-1}[v_1xyx^{\q-2}]]
\\&=
\binom{\q-2}{q-2}[v_{-1}[v_1xy]x^{q}]
+\binom{\q-2}{q-1}[v_{-1}x[v_1xy]x^{q-1}]
\\&=
[v_{-1}[v_1x]yx^{q}]
\\&=
[v_0yx^{q}].
\end{align*}
We conclude that $[v_0yx^{q}]=0$,
and hence $\dim L\le k+q+1$.

We now may assume that the sequence of constituent lengths of $M$ reads
\[\boxed{
\q,\q-1, \q^{r},\ldots,\q-1,\q^{r-2},\q-1,
}\]
and hence $[v_2[v_1x]^rx]\neq 0$, but $[v_2[v_1x]^ry]=0$.
We continue to assume as we may that $r$ is even, and that the penultimate constituent of $M$ of length $\q-1$
is immediately preceded by a constituent of length $\q$.
Accordingly, we redefine $v_{-r-1}$ to be
the unique homogeneous element such that
$[v_{-r-1}xyx^{\q-2}]=v_r$.
In order to bound $\dim L$ we need to make a further case distinction, according as the
penultimate constituent of $M$ of length $\q-1$ is immediately preceded by the pair $\q-1,\q$ or by $\q,\q$,
and we define an element $v_{-r-2}$ accordingly.

In the former case, the sequence of constituent lengths of $M$ reads
\[\boxed{
\q,\q-1, \q^{r},\ldots,\q-1,\q,\q-1,\q^{r-2},\q-1,
}\]
and we define $v_{-r-2}$ to be
the unique homogeneous element such that
$[v_{-r-2}xyx^{\q-3}]=v_{-r-1}$.
Because $r\ge 2$ we have $[v_2[v_1x]^2y]=0$, and hence
\begin{align*}
0&=
[v_{-r-2}x[v_2[v_1x]^2y]]
\\&=
[v_{-r-2}x[v_2[v_1x]^2]y]
+[v_{-r-2}xy[v_2[v_1x]^2]]
\\&=
[v_{-r-2}xv_2[v_1x]^2y]
+[v_{-r-2}xy[v_1x]^2v_2]
\\&=
[v_{-r}[v_1x]^2y]
+[v_{-r}xyxv_2]
\\&=
[v_{-r+2}xy]
+[v_{-r}xyxv_2],
\end{align*}
where we have used the facts that
$[v_{-r-2}xv_2]=v_{-r}$
(which is analogous to Equation~\eqref{eq:v_-rxv_2}), and
$[v_{-r-2}xyv_2]=0$.
When $r>2$ we have
\[
[v_{-r}xyxv_2]
=
[v_{-r}xyx^{\q-2}[v_1xy]]
=
[v_{-r+1}x[v_1x]y]+
[v_{-r+1}xy[v_1x]]
=0,
\]
and when $r=2$ we have
\begin{align*}
[v_{-r}xyxv_2]
&=
[v_{-1}[v_1xy]x]
+[v_{-1}x[v_1xy]]
\\&=
[v_{-1}[v_1x]yx]
+[v_{-1}x[v_1x]y]
+[v_{-1}xy[v_1x]]
\\&=
[v_0yx]+[v_0xy]+\mu^{-1}[v_0xy]
=0.
\end{align*}
In all cases we conclude that $[v_{-r+2}xy]=0$, which is a contradiction.

It remains to deal with the case where the sequence of constituent lengths of $M$ reads
\[\boxed{
\q,\q-1, \q^{r},\ldots,\q-1,\q^{r-2},\q-1,
}\]
and the penultimate constituent of length $\q-1$ is immediately preceded by at least two constituents of length $\q$.
Here we redefine $v_{-r-2}$ to be
the unique homogeneous element such that
$[v_{-r-2}xyx^{\q-2}]=v_{-r-1}$.
Therefore, differently from the case of the previous paragraph, here we have
$[v_{-r-2}xv_2]=0$
(as in Equation~\eqref{eq:v_mxv_2-bis}).
We now employ the fact that
$[v_2[v_1x]^ry]=0$.
Recalling that $r$ is even and assuming $r>2$ we expand
\begin{align*}
0&=
[v_{-r-2}x[v_2[v_1x]^ry]]
\\&=
[v_{-r-2}x[v_2[v_1x]^r]y]
+[v_{-r-2}xy[v_2[v_1x]^r]]
\\&=
\binom{r}{2}
[v_{-r-2}x[v_1x]^2v_2[v_1x]^{r-2}y]
+[v_{-r-2}xy[v_1x]^rv_2]
\\&=
\binom{r}{2}
[v_{-r}xv_2[v_1x]^{r-2}y]
+[v_{-2}xyxv_2]
\\&=
\binom{r}{2}[v_0y]
+[v_0y].
\end{align*}
In this calculation we have also used the facts that
$[v_{-r}xv_2]=v_{-r+2}$
(which is Equation~\eqref{eq:v_-rxv_2}),
$[v_{-r-2}xyv_2]=0$,
$[v_{-r-1}xyv_2]=0$,
$[v_{-r}xyv_2]=0$,
$[v_mxyxv_2]=0$ for $-r<m<-2$
and
$[v_{-2}xyxv_2]=[v_0y]$
(which is Equation~\eqref{eq:v_-1xyxv_2}, valid for $r>2$).
Note that these formulas require no calculation except for the first and the last,
if we recall that $L$ is $\Z^2$-graded by $x$ and $y$.
If $r$ is a multiple of four, then $r>2$ and $\binom{r}{2}$ is even.
Hence the above calculation applies and yields a contradiction.

It remains to deal with the case where $r$ (is even but) is not a multiple of four, so that $\binom{r}{2}$ is odd.
Assuming $r>2$ for now, we expand
\begin{equation}\label{eq:final}
\begin{split}
0&=
[v_{-r}x[v_2[v_1x]^ry]]
\\&=
[v_{-r}x[v_2[v_1x]^r]y]
+[v_{-r}xy[v_2[v_1x]^r]]
\\&=
[v_{-r}xv_2[v_1x]^ry]
+\binom{r}{2}
[v_{-r}xy[v_1x]^{r-2}v_2[v_1x]^2]
+[v_{-r}xy[v_1x]^{r}v_2]
\\&=
[v_{-r+2}[v_1x]^ry]
+\binom{r}{2}
[v_{-2}xyxv_2[v_1x]^2]
+[v_{-1}xyx[v_1x]v_2]
\\&=
[v_0[v_1x]^2y]
+\binom{r}{2}
[v_0y[v_1x]^2]
+\mu^{-1}[v_0xyxv_2]
\\&=
\mu^{-1}[v_{n+2}xy]
+\binom{r}{2}
[v_{n+2}xy],
\end{split}
\end{equation}
because
$[v_0xyxv_2]=[v_{n+1}[v_1xy]]=0$
and other formulas already used.
Therefore, when $r>2$ is not a multiple of four
we obtain $(\mu^{-1}+1)[v_{n+2}xy]=0$, whence $[v_{n+2}xy]=0$.
Because $[v_{n+2}xx]=0$ as well, we conclude that $\dim L\le k+2\q+1$.
When $r=2$, Equation~\eqref{eq:final} needs to be modified, but before doing that
we show that
$[v_{n+2}y]=0$
in this case as well
(which only follows from our adaptation of Lemma~\ref{lemma:8} when $r>2$).
In fact, because
\begin{align*}
[v_{-2}v_2]
&=
[v_{-2}[v_1xy]x^{\q-3}]
+\binom{\q-3}{1}[v_{-2}x[v_1xy]x^{\q-4}]
\\&=
[v_{-2}[v_1x]yx^{\q-3}]
+[v_{-2}xy[v_1x]x^{\q-4}]
\\&=
[v_{-1}xyx^{\q-3}]
+[v_{-1}xyx^{\q-3}]
=0,
\end{align*}
we have
\[
0=[v_{-2}[v_2[v_1x]^2y]]
=[v_{-2}[v_1x]^2v_2y]
=[v_0xv_2y]
=[v_{n+2}y]
\]
according to Equation~\eqref{eq:v_nxv_2}.
Now the analogue of Equation~\eqref{eq:final} is
\begin{align*}
0&=
[v_{-r}x[v_2[v_1x]^ry]]
\\&=
[v_{-2}x[v_2[v_1x]^2]y]
+[v_{-2}xy[v_2[v_1x]^2]]
\\&=
[v_{-2}xv_2[v_1x]^2y]
+[v_{-2}xyv_2[v_1x]^2]
+[v_{-2}xy[v_1x]^{2}v_2]
\\&=
[v_0[v_1x]^2y]
+[v_{-1}[v_1xy][v_1x]^2]
+[v_{-1}xyx[v_1x]v_2]
\\&=
[v_0[v_1x]^2y]
+[v_0y[v_1x]^2]
+\mu^{-1}[v_0xyxv_2]
\\&=
\mu^{-1}[v_{n+2}xy]
+(\mu^{-1}+1)[v_{n+2}xy]
\\&=
[v_{n+2}xy]
\end{align*}
and implies, again, that $\dim L\le k+2\q+1$.

Noting that none of the upper bounds for $\dim L$ which we have found, in the various cases, exceeds
$(4k+1)/3$, completes the proof of Theorem~\ref{thm:main}.

\bibliography{References}

\def\cprime{$'$} \def\cprime{$'$}
\providecommand{\bysame}{\leavevmode\hbox to3em{\hrulefill}\thinspace}
\providecommand{\MR}{\relax\ifhmode\unskip\space\fi MR }
\providecommand{\MRhref}[2]{%
  \href{http://www.ams.org/mathscinet-getitem?mr=#1}{#2}
}
\providecommand{\href}[2]{#2}
\begin{thebibliography}{CMNS96}

\bibitem[AJ01]{AviJur}
M.~Avitabile and G.~Jurman, \emph{Diamonds in thin {L}ie algebras}, Boll.
  Unione Mat. Ital. Sez. B Artic. Ric. Mat. (8) \textbf{4} (2001), no.~3,
  597--608. \MR{MR1859998 (2003a:17038)}

\bibitem[AM05]{AviMat:-1}
M.~Avitabile and S.~Mattarei, \emph{Thin {L}ie algebras with diamonds of finite
  and infinite type}, J. Algebra \textbf{293} (2005), no.~1, 34--64.
  \MR{MR2173965 (2006f:17018)}

\bibitem[AM07]{AviMat:A-Z}
\bysame, \emph{Thin loop algebras of {A}lbert-{Z}assenhaus algebras}, J.
  Algebra \textbf{315} (2007), no.~2, 824--851. \MR{MR2351896 (2008h:17022)}

\bibitem[Avi02]{Avi}
M.~Avitabile, \emph{Some loop algebras of {H}amiltonian {L}ie algebras},
  Internat. J. Algebra Comput. \textbf{12} (2002), no.~4, 535--567.
  \MR{MR1919687 (2003e:17013)}

\bibitem[Bla58]{Blackburn}
N.~Blackburn, \emph{On a special class of {$p$}-groups}, Acta Math.
  \textbf{100} (1958), 45--92. \MR{MR0102558 (21 \#1349)}

\bibitem[Bra88]{Br}
Rolf Brandl, \emph{The {D}ilworth number of subgroup lattices}, Arch. Math.
  (Basel) \textbf{50} (1988), no.~6, 502--510. \MR{MR948264 (89e:20054)}

\bibitem[Cam00]{Cam}
Rachel Camina, \emph{The {N}ottingham group}, New horizons in pro-$p$ groups,
  Progr. Math., vol. 184, Birkh\"auser Boston, Boston, MA, 2000, pp.~205--221.
  \MR{MR1765121 (2001f:20054)}

\bibitem[Car97]{Car:Nottingham}
A.~Caranti, \emph{Presenting the graded {L}ie algebra associated to the
  {N}ottingham group}, J. Algebra \textbf{198} (1997), no.~1, 266--289.
  \MR{MR1482983 (99b:17019)}

\bibitem[Car99]{Car:Zassenhaus-three}
\bysame, \emph{Loop algebras of {Z}assenhaus algebras in characteristic three},
  Israel J. Math. \textbf{110} (1999), 61--73. \MR{MR1750443 (2001d:17019)}

\bibitem[CJ99]{CaJu:quotients}
A.~Caranti and G.~Jurman, \emph{Quotients of maximal class of thin {L}ie
  algebras. {T}he odd characteristic case}, Comm. Algebra \textbf{27} (1999),
  no.~12, 5741--5748. \MR{MR1726275 (2001a:17042a)}

\bibitem[CM99]{CaMa:thin}
A.~Caranti and S.~Mattarei, \emph{Some thin {L}ie algebras related to
  {A}lbert-{F}rank algebras and algebras of maximal class}, J. Austral. Math.
  Soc. Ser. A \textbf{67} (1999), no.~2, 157--184, Group theory. \MR{MR1717411
  (2000j:17036)}

\bibitem[CM04]{CaMa:Nottingham}
\bysame, \emph{Nottingham {L}ie algebras with diamonds of finite type},
  Internat. J. Algebra Comput. \textbf{14} (2004), no.~1, 35--67. \MR{MR2041537
  (2004j:17027)}

\bibitem[CM05]{CaMa:Hamiltonian}
\bysame, \emph{Gradings of non-graded {H}amiltonian {L}ie algebras}, J. Aust.
  Math. Soc. \textbf{79} (2005), no.~3, 399--440. \MR{MR2190690 (2006j:17020)}

\bibitem[CMN97]{CMN}
A.~Caranti, S.~Mattarei, and M.~F. Newman, \emph{Graded {L}ie algebras of
  maximal class}, Trans. Amer. Math. Soc. \textbf{349} (1997), no.~10,
  4021--4051. \MR{MR1443190 (98a:17027)}

\bibitem[CMNS96]{CMNS}
A.~Caranti, S.~Mattarei, M.~F. Newman, and C.~M. Scoppola, \emph{Thin groups of
  prime-power order and thin {L}ie algebras}, Quart. J. Math. Oxford Ser. (2)
  \textbf{47} (1996), no.~187, 279--296. \MR{MR1412556 (97h:20036)}

\bibitem[CN00]{CN}
A.~Caranti and M.~F. Newman, \emph{Graded {L}ie algebras of maximal class.
  {II}}, J. Algebra \textbf{229} (2000), no.~2, 750--784. \MR{MR1769297
  (2001g:17041)}

\bibitem[Ers05]{Ershov:finitely_presented}
M.~V. Ershov, \emph{The {N}ottingham group is finitely presented}, J. London
  Math. Soc. (2) \textbf{71} (2005), no.~2, 362--378. \MR{MR2122434
  (2006c:20058)}

\bibitem[GAP07]{GAP}
The GAP Group, \emph{\textsf{GAP} --- {G}roups, {A}lgorithms, and
  {P}rogramming, version 4.4.10}, 2007, \texttt{http://www.gap-system.org}.

\bibitem[GMY01]{GMY}
Norberto Gavioli, Valerio Monti, and David~S. Young, \emph{Metabelian thin
  {L}ie algebras}, J. Algebra \textbf{241} (2001), no.~1, 102--117.
  \MR{MR1838846 (2002d:17032)}

\bibitem[Hup67]{Hup}
B.~Huppert, \emph{Endliche {G}ruppen. {I}}, Die Grundlehren der Mathematischen
  Wissenschaften, Band 134, Springer-Verlag, Berlin, 1967. \MR{MR0224703 (37
  \#302)}

\bibitem[Jur99]{Ju:quotients}
G.~Jurman, \emph{Quotients of maximal class of thin {L}ie algebras. {T}he case
  of characteristic two}, Comm. Algebra \textbf{27} (1999), no.~12, 5749--5789.
  \MR{MR1726276 (2001a:17042b)}

\bibitem[Jur04]{Ju:Bi-Zassenhaus}
\bysame, \emph{A family of simple {L}ie algebras in characteristic two}, J.
  Algebra \textbf{271} (2004), no.~2, 454--481. \MR{MR2025538 (2005a:17015)}

\bibitem[Jur05]{Ju:maximal}
\bysame, \emph{Graded {L}ie algebras of maximal class. {III}}, J. Algebra
  \textbf{284} (2005), no.~2, 435--461. \MR{MR2114564 (2005k:17041)}

\bibitem[JY]{JuYo:quotients}
G.~Jurman and D.~S. Young, \emph{Quotients of maximal class of thin {Lie}
  algebras in characteristic two: errata and addendum}, unpublished manuscript,
  2003.

\bibitem[KLGP97]{KL-GP}
G.~Klaas, C.~R. Leedham-Green, and W.~Plesken, \emph{Linear pro-{$p$}-groups of
  finite width}, Lecture Notes in Mathematics, vol. 1674, Springer-Verlag,
  Berlin, 1997. \MR{MR1483894 (98m:20028)}

\bibitem[LGM02]{L-GMcKay}
C.~R. Leedham-Green and S.~McKay, \emph{The structure of groups of prime power
  order}, London Mathematical Society Monographs. New Series, vol.~27, Oxford
  University Press, Oxford, 2002, Oxford Science Publications. \MR{MR1918951
  (2003f:20028)}

\bibitem[Mat]{Mat:chain_lengths}
S.~Mattarei, \emph{Constituents of graded {L}ie algebras of maximal class and
  chain lengths of thin {L}ie algebras}, in preparation.

\bibitem[Mat99]{Mat:thin-groups}
S.~Mattarei, \emph{Some thin pro-{$p$}-groups}, J. Algebra \textbf{220} (1999),
  no.~1, 56--72. \MR{MR1713453 (2000h:20049)}

\bibitem[You01]{Young:thesis}
D.~S. Young, \emph{Thin {Lie} algebras with long second chains}, Ph.D. thesis,
  Canberra, March 2001.

\end{thebibliography}

\end{document}